\documentclass[reqno]{amsart}
\usepackage{amssymb}
\usepackage{epsf}
\usepackage{xcolor}
\usepackage{hyperref}

\theoremstyle{definition}

\theoremstyle{remark}

\numberwithin{equation}{section}

\newcommand{\be}{\begin{equation}}
\newcommand{\ee}{\end{equation}}
\newcommand{\ba}{\begin{eqnarray}}
\newcommand{\ea}{\end{eqnarray}}
\newcommand{\baa}{\begin{eqnarray*}}
\newcommand{\eaa}{\end{eqnarray*}}
\newcommand{\bb}{\begin{thebibliography}}
\newcommand{\eb}{\end{thebibliography}}

\newcommand{\bi}[1]{\bibitem{#1}}
\newcommand{\lab}[1]{\label{#1}}
\newcommand{\re}[1]{(\ref{#1})}

\newcommand{\ve}{\varepsilon}

\newtheorem{tm}{Theorem}
\newtheorem{pr}{Proposition}
\newtheorem{cor}{Corollary}
\newtheorem{lem}{Lemma}
\newtheorem{st}{Statement}
\renewcommand{\t}{\tilde}

\newcommand{\Q}{\mathbb Q}
\newcommand{\C}{\mathbb C}
\newcommand{\R}{\mathbb R}
\newcommand{\Z}{\mathbb Z}

\begin{document}


\title[ On Dirichlet, Poncelet and Abel problems] 
{On Dirichlet, Poncelet and Abel problems}

\author{V.P. Burskii }
\address{Institute of Applied Mathematics, Donetsk, 83114, Ukraine}

\author{A.S. Zhedanov}
\address{Donetsk Institute for Physics and Technology,
Donetsk, 83114, Ukraine}


\begin{abstract}
We propose interconnections between some problems of PDE, geometry, algebra, calculus and physics. Uniqueness of a solution of the Dirichlet problem and of some other boundary value problems for the string equation inside an arbitrary biquadratic algebraic curve is considered. It is shown that a solution is non-unique if and only if a corresponding Poncelet problem for two conics has a periodic trajectory. A set of problems is proven to be equivalent to the above problem. Among them are the solvability problem of the algebraic Pell-Abel equation and the indeterminacy problem of a new moment problem that generalizes the well-known trigonometrical moment problem. Solvability criteria of the above-mentioned problems can be represented in form $\theta\in\Bbb Q$ where number $\theta=m/n$ is built by means of data for a problem to solve. We also demonstrate close relations of the above-mentioned problems to such problems of modern mathematical physics as elliptic solutions of the Toda chain, static solutions of the classical Heisenberg $XY$-chain and biorthogonal rational functions on elliptic grids in the theory of the Pad\'e interpolation.
\\[2mm]
\end{abstract}

\maketitle

\tableofcontents


\section{Introduction}
\setcounter{equation}{0}

This work is dedicated to establishing some new interconnections between certain problems of mathematics such as ill-posed boundary value problems in bounded semi-algebraic domains for partial differential equations, the moment problem, the Poncelet problem from projective geometry, the algebraic Pell-Abel equation and some other problems, recently revealed by the authors.

Study of ill-posed boundary value problems in bounded domains for partial differential equations goes back to J. Hadamard. These problems serve as a popular object of present-day investigations (see s. \ref{HisRemBVP}). In this paper we will examine a general Dirichlet problem for the string equation
\be u_{xy} =0 \quad
\mbox{in}\ \ \Omega, \lab{SE} \ee \be u|_{C} =\phi \quad
\mbox{on}\ \ C=\partial \Omega,\lab{NDP} \ee
for which the solution uniqueness is related to a problem of existence of a nontrivial solution of the homogeneous Dirichlet problem.
\be
u|_{C} =0 \quad \mbox{on}\ \ C.\lab{NDPH} \ee
Functions
$u,\,\phi$
are assumed to be complex-valued functions of real variables. We will consider this problem as well as some other boundary value problems in semi-algebraic domains, the boundaries of which are given by so-called biquadratic algebraic curves
\be F(x,y):=\sum_{i,k=0}^2
a_{ik}x^iy^k=0. \lab{2_2F} \ee
We will consider canonical forms of curve \re{2_2F}, to which the generic curve can be transformed by linear-fractional replacements, and we will come up with criteria of uniqueness breakdown in the form
\be\label{crit}\tau\in\Bbb Q\ee
where number $\tau$ is determined by curve $C$. Our investigations are based on a study of the John mapping generated on $C$ by characteristics of equation \re{SE}; see proposition \ref{CondBreak} (s. \ref{breakdown}) and proposition \ref{etaJohn} (s. \ref{Jbc}) in which we use that and which reflect that the John mapping becomes an ordinary shift after a transform on universal covering group of variety \re{2_2F}.

We discovered a remarkable similarity of this problem with the Poncelet problem. The Poncelet problem is one of well-known problems of projective geometry and the problem on its own possesses numerous links to a set of various problems of calculus, functional analysis and physics (see \cite{Berger},\cite{Bax} \cite{GZ} and s. \ref{HisRemPP} below). We will prove that for generic biquadratic curve $C$ the Dirichlet problem has a non-unique solution if and only if a corresponding Poncelet problem has a periodic trajectory (proposition \ref{PoncJohn}, s. \ref{PerPonc}). And therefore we will give a new criterion for periodicity in the form \re{crit} that differs from the well-known Cayley criterion (see s. \ref{Cayley}).

To remind, the existence of a periodic trajectory in the Poncelet problem implies that every trajectory is periodic in accordance with Poncelet big theorem \cite{Berger}. Note that different cases of disposition of conics give different cases of curves $C$ and a new setting of the Dirichlet problem \re{SE},\re{NDPH} for unbounded domains is more suitable than the classical setting (see s. \ref{JohnC}).

We will observe one more remarkable similarity of the Poncelet problem with the solvability problem of the algebraic Pell-Abel equation
\be P(t)^2-R(t)Q(t)^2=L,\lab{PellAb}\ee
where for a given polynomial $R(t)$ of one variable $t$ one seeks polynomials $P, Q$ and a constant $L$ satisfying the equation. Solvability of the Pell-Abel equation is also one of well-known problems and this algebraic problem by itself has numerous connections to many problems of calculus, functional analysis and function theory (see s. \ref{PAbel} below).
We will prove (see proposition \ref{t_PA_Poncelet} of s. \ref{PAPoncelet}) that the Pell-Abel equation \re{Pell} with a polynomial $R$ of the fourth order has a solution if and only if a corresponding Poncelet problem has a periodic trajectory of an even period.
Doing so, we obtain a new criterion in the form \re{crit} which differs from already existing well-known criteria such as the well-known Zolotarev porcupine.

We will show that the same condition \re{crit} serves as a criterion of solution uniqueness to within an additive constant of the Neumann problem
\be u_{\nu_*}|_{C} =\psi \lab{NNP} \ee
for the same equation \re {SE} in the same domain where $u_{\nu_*}$ is a derivative with respect to the conormal $\nu_*$ (statement \ref{DirNeuMom} of s. \ref{DirNeuMo}). We will also show that the same condition serves as a criterion of indeterminacy of a new moment problem on curve $C$ (statement \ref{DirNeuMom} of s. \ref{DirNeuMo}).
Through the equation-domain duality we will obtain an equivalent problem \re{DBEEE} in the form of a hyperbolic equation (in certain cases of a corresponding Poncelet problem) of the fourth order with only two boundary data on characteristics instead of four boundary data, as it would be in the boundary value problem of the Goursat type, and like the Goursat type boundary value problems, the solution for our problem is also almost always unique.

We will observe some similarities with a problem of classical $XY$-spin chains and a problem from the theory of the Toda chain. We will show an equivalence of the considered problems and also provide an interpretation of this criterion in terms of the John mapping. For the links to some other problems of mathematics, see \cite{Ves}, \cite {Laz}, \cite {Malyshev}, \cite{SodinYud}.

Note that some results of the present work were already published in the form of brief and incomplete fragments in papers \cite{BurZhedan1},
\cite{BurZhedan2} and \cite{BurZhedan3}.

Note also that some explicit necessary and sufficient conditions of uniqueness breakdown of solution of the Dirichlet problem (and some others boundary value problems) for partial differential equations with constant coefficients were obtained earlier for an arbitrary ellipse (see e.g., \cite{BurMomProblem} and \cite{BurBook}). Answers in those works were formulated in the form of condition (1.5), which was to be a hint in our present investigations. In this paper we would like to offer a new way of examining ill-posed boundary value problems for partial differential equations in domains, which are more complex than a circle, along with surveying some other fields of mathematics for topics of equivalent contents.

Hence, the primary objective of the present paper is to show numerous similarities and relations of the considered ill-posed Dirichlet problem with remotely located branches of mathematics and mathematical physics.

The paper is organized as follows. In Section 2, a theory of uniqueness of the Dirichlet problem and the John algorithm are described. In section 3, application of this technique to certain biquadratic curves is provided. In Section 4, solution of the Poncelet problem is considered and its relation to the John algorithm for a generic biquadratic curve is shown. Section 5 is dedicated to the Abel problem of explicit computation of elliptic integrals by means of elementary functions. It is shown how the problem can be formulated in terms of the Poncelet problem (or of the John algorithm for a biquadratic curve). In Section 6, a relation of our problems to Ritt's problem of existence of periodic functions with a nontrivial multiplication property is considered. Finally, in Section 7, we consider three problems of mathematical physics which are related to the Poncelet problem (or the John algorithm): static solutions of the classical $XY$ Heisenberg chain, elliptic solutions of the Toda chain and elliptic grids for biorthogonal rational functions in the theory of rational Pad\'e interpolation.

\vskip 25pt

\section{Boundary value problems in
a domain for the string equation}

\subsection{
Bibliographical remarks.}\label{HisRemBVP}
Investigations of ill-posed boundary value problems in bounded domains for partial differential equations go back to J. Hadamard \cite{Had} and further to A. Huber \cite{Hub}, who were the first to note non-uniqueness of the solution of the Dirichlet problem for the equation of string vibration (i.e. the string equation) in the rectangle. Boundary value problems in bounded domains for non-elliptic partial differential equations were studied regularly, for the most part, for a parallelepiped. Likewise, boundary value problems were also studied more or less regularly in domains with general boundary questions of solution uniqueness of the Dirichlet problem for a hyperbolic equation of the second order on a plane (see reviews and results in \cite{Pta}). In work \cite{BD}, D. Burgin and R. Duffin examined a homogeneous Dirichlet problem for equation $u _{tt}-u _ {xx} =0 $ in rectangular $ \{0\le t\le T; 0\le x\le X
\} $. It was shown that if the ratio $T/X$ is irrational, uniqueness in space of continuously differentiated functions with summable second derivatives takes place. Theorems of existence of solutions in classical spaces were developed, and it was shown that the greater the smoothness of a boundary function was, and likewise, the worse the number $T/X$ was approximated by means of rational numbers, the greater the smoothness of the solution was observed to take place, as a result. The Neumann problem was considered there also. In works of B. Yo. Ptashnik and his disciples, boundary value problems inside a parallelepiped for a wide class of differential equations and systems of differential equations were investigated (see \cite{Pta}). All those works (except for the mentioned work \cite{Arnold}) were based on methods which essentially use the representation of such domain in view of the topological product.

For nonrectangular domains, the Dirichlet problem for the string equation was studied in connection with the number of the Denjoy-Poincare rotation (see, e.g., Z. Nitetsky's book \cite{Nit}) of homeomorphism of the domain boundary, constructed on characteristics of the equation (so-called automorphism of characteristic billiards by Fritz John \cite{John}). A connection of properties of the Dirichlet problem to properties of this homeomorphism was used even in the above mentioned works of J. Hadamard and A. Huber. When it comes to calculus, this connection was established in F. John's work \cite{John}. In works by R. A. Alexandrjan and his disciples, investigations of this problem and, in particular, of this connection were continued (\cite{Alex1},\cite{Alex2}, \cite{Ovs} and
\cite{Alex3}). The question of uniqueness of the solution of the Dirichlet problem in this ideology for domains that are convex with respect to characteristics families should be transformed to the question of irrationality of rotation number or, identically, to the question of presence a continuous set of finite orbits (cycles) of the discrete dynamic system generated by the above-mentioned John homeomorphism.

The same questions in connection with an asymptotic behavior of the solution of the Sobolev equation that describe surface oscillations of a fluid filling a body, flying through the atmosphere, were investigated by Siberian mathematicians T. I. Zelenjak, I. V. Fokin and some others (\cite{Zell}, \cite{Fokin}). Research studies of the string equation were also included into the well-known book of Yu. M. Berezansky \cite{Ber}; the results from the studies describe possibilities for building domains with angles, at which the homogeneous Dirichlet problem is weakly solvable and well-posed, concerning the right part and small motions of a boundary of the domain, with angles left within specified limits. Note also that the case of an ellipse was considered in works by A. Huber  \cite{Hub}, R. Alexandrjan \cite{Alex2}, V. I. Arnold
\cite{Arnold}. For information on small smoothness and on more general equations, see the book \cite{BurBook}.

\subsection[John]{John condition}\label{JohnC}

For problem \re{SE}, \re{NDP}, solved for some general bounded domains, Fritz John \cite{John} considered a remarkable transformation $T:C\to C$ of the Jordan boundary into itself, allowing for making some conclusions regarding properties of the Dirichlet problem in $\Omega$. Let us describe it in more detail.

Let $\Omega$ be an arbitrary bounded domain which is convex with respect to characteristics of the equation \re{SE}, i.e. it has the boundary $C$ intersected, at most, in two points by every straight line which is parallel to the $x$- or $y$-axes. We start from arbitrary point $M_1$ on $C$ and consider a vertical line passing through $M_1$. Generally, there are two intersection points with curve $C$: $M_1$ and some point $M_2$. We denote as $I_1$ involution which transforms $M_1$ into $M_2$. Then, starting from $M_2$, we consider a horizontal line passing through $M_2$. Let $M_3$ be the second point of intersection with curve $C$. Let $I_2$ be corresponding involution:
$I_2M_2=M_3$.

Then we repeat this process, applying, step-by-step, involutions $I_1$ and $I_2$. Denote $T=I_2I_1, T^{-1}=I_1I_2$. This transformation $T:C\to C$ gives a discrete dynamical system on $C$, i.e. an action of group $\Z$ and each point $M\in C$ generates an orbit $\{T^nM|n\in\Z \}$. This orbit can be a finite or denumerable set. The point $M$ with a finite orbit is called a periodic point and the smallest $n$, for which $T^nM=M$, is called a period of point $M$. In the paper \cite{John}, uniqueness breakdown in problems \re{SE},\re{NDP} were studied in connection with topological properties of the mapping $T$ for the case of the even mapping $T$. The mapping $T$ is called even or preserving an orientation, if every positively oriented arc $(P,Q)$ with points $P,Q\in C$ transforms into a positively oriented arc $(TP,TQ)$. F. John has proven several useful assertions, from which we extract the following.

{\it Sufficient condition of uniqueness.} The homogeneous Dirichlet problem \re{SE}, \re{NDPH} in a bounded domain has only a trivial solution in space $C^2(\overline\Omega)$ if a set of periodic points on $C$ is finite or denumerable.

Four possible cases of dynamical system behavior are selected here, which are

I) all points are periodic (in which case, their periods coincide);

II) there exist periodic and nonperiodic points;

III) there are no periodic points and there exists no point, an orbit of which $\{...,T^{-1}P,P,TP,...,T^nP,...\}$ is dense on $C$;

IV) there are no periodic points and there exists a point, an orbit of which $\{...,T^{-1}P,P,TP,...,T^nP,...\}$ is dense on $C$ (a transitive case).

\vskip 5pt
In the work \cite{John}, it was shown that for a $C^2$-smooth curve, case III may not be present. For case II, it was proven that there exists an arc $D_0$ on $C$ such that every two arcs from $D_0,I_1D_0,TD_0,I_1TD_0,T^2D_0,$ ..., $I_1T^{n-1}D_0,T^nD_0,...$ have no common points. Note that for an analytical boundary, this may not be valid, because in this case $T$ is a diffeomorphism. In case IV, it was proven that the Denjoy-Poincar\'e rotation $\xi$ (\cite{Nit}) of the John mapping $T$ is irrational and $T$ is topologically equivalent to a turn of a unit circle about the angle $\pi\xi$ (i.e. there exists a homeomorphism $h$ from $C$ onto a unit circle $S$ such that the mapping $hTh^{-1}:S\to S$ is a turn about the angle $\pi\xi$).
For a case, in which every point of $C$ is periodic, it was proven that all periods coincide. But in this case, nothing is known about solution uniqueness, although R.A. Alexandrjan has shown for this case \cite{Alex1} that there is a generalized solution of the problem \re {SE},\re{NDP}, which can be constructed in the form of a piecewise-constant function.

We will assume that the domain has boundary $C$ satisfying the condition:
\begin{eqnarray}\label{CoDom}\mbox {\it The curve $C$ is smooth and each characteristic
line either doesn't intersect\qquad}\nonumber \\ \mbox{\it
the curve $C$ or 
touches it at a point, called a vertex, or splits $C$ in two
points.\ \ 
}\end{eqnarray}

Note that John's condition of convexity with respect to characteristic directions given above will be satisfied under the condition \re{CoDom} on $C=\partial\Omega$ in the case of a bounded domain
$\Omega$.

For the case of a bounded domain with a biquadratic boundary, we will see that the sufficient uniqueness condition given above is also necessary, and moreover, it will be so even for cases, in which curve $C$ is unbounded, but in that case we should change the setting of the problem. Namely, along with the usual setting of the uniqueness property:
\begin{eqnarray}\label{UHDP} \mbox{ \it an examined bounded domain
is such that the homogeneous Dirichlet \ \ \ \qquad}\nonumber \\
\mbox{\it \ problem \re{SE},\re{NDPH} has only a trivial solution in
the space $C^2(\Omega)\cap C(\overline\Omega)$,\
\qquad}\end{eqnarray} for cases, in which the curve $C$ is unbounded, we
will examine the following modification of the uniqueness property
for the homogeneous Dirichlet problem
:
\begin{eqnarray}\label{MHDP} \mbox{ \it  an examined curve $C$ is
such that every analytic, in real sense solution \ \ \ \ }\nonumber \\
\mbox{ \it in $\Bbb R^2$ of the equation \re{SE} with the property
\re{NDPH} may be only a zero solution.\quad }\end{eqnarray}

Note that here the assumption of an "analytic" curve is introduced so that now we can consider unbounded curves $C$, on which there exist characteristic lines which do not intersect $C$ and which are located between curve branches. Without this assumption on solutions for such a curve and, for example, using an assumption of infinite smoothness for solutions, one may build a simple example of a smooth in $\Bbb R^2$ nontrivial solution of equation \re{SE} with the property \re{NDPH}. We will call a modified setting of the problem to solve as an analytic in the real sense solution in $\Bbb R^2$ of equation \re{SE} with property \re{NDPH}.

Now we will give John's proof of the sufficient condition of uniqueness in order to show that John's arguments are also valid for the problem \ref{MHDP}.

\begin{pr}\label{unique}
Under the condition \re{CoDom}, if the mapping $T$ transitively acts
on $C$, then 1) the uniqueness property \re{UHDP} holds for an
ordinary setting of the homogeneous Dirichlet problem
\re{SE},\re{NDPH} in a bounded domain and 2) the uniqueness
property \re{MHDP} holds for the modified setting in a, possible,
unbounded domain.
\end{pr}

\begin{proof}
Let a function $u\in C^2(\Omega)\cap C(\overline\Omega)$ be a nontrivial solution of the problem \re{SE},\re{NDPH} in domain $\Omega$ with the condition \re{CoDom}. As is well-known, there exist two functions $u_1,u_2$ of the class $C^2$ depending on one variable such that $u(x,y)=u_1(x)+u_2(y)$, which we will write down for any point $ P\in C$ as $u(P)=u_1(P)+u_2(P)$. For the case of the property \re{MHDP}, consider a domain $\Omega\subset \Bbb R^2$ of points $P=(x_0,y_0)\in\Bbb R^2$ for which there exists a pair of different points of intersection $C\cap \{x=x_0\}$ of the curve with a corresponding vertical line and also there exists a pair of different points of intersection $C\cap \{y=y_0\}$ of the curve with a horizontal characteristic line.
Then for any point $P\in\Omega$, using the definitions, we can easily obtain:
$u_1(I_1P)=u_1(P),$ $u_2(I_2P)=u_2(P);$
$u_2(P)-u_2(TP)=u_2(P)-u_2(I_1P)=u(P)-u(I_1P)=0;$
$u_1(P)-u_1(T^{-1}P)=u_1(P)-u_1(I_2P)=u(P)-u(I_2P)=0.$

From that it follows that equalities $u_1(P)=u_1(T^nP)$, $u_2(P)=u_2(T^nP)$ hold for every integer $n$. Continuity gives us $u_1(P)=u_1(Q)$, $u_2(P)=u_2(Q)$ for any point $Q$ from closure of the orbit of $P$. Because of a transitive acting, the closure of the orbit of any point coincides with $C$, then $u_1\equiv const$ and $u_2\equiv const$ and, therefore, $u\equiv 0$ in
$\overline\Omega$. In the case of the property \re{MHDP}, the analyticity allows us to extend the zero solution onto the plane.

\end{proof} 

It was noted above that for any analytic boundary, only two cases are admissible: periodic (I) and transitive (IV by John). Therefore we give the following setting of {\it the periodicity problem for John mapping} :

What curve from a given class of curves has the following property?
\be\label{JPerCond}\quad The \ John\ mapping\ \ T:C\to C\ \ has\
at\ least\ one\ periodic\ point.\ee
Then, as we know already, all the points on our curve are periodic.

Along with the above settings, we will consider also a setting with the complex John mapping. For a biquadratic complex curve (i.e. one-dimensional complex variety \re{2_2F}) $\tilde C\in \Bbb C^2$, we will assume that it is a satisfied property of type \ref{CoDom}:
\begin{eqnarray*}\label{Complconv} 
\mbox{\it Almost every "vertical"
line } x=x_0\ \mbox{\it intersects the curve \ }
\tilde C \mbox{\it \ at two different points}\nonumber \\
\mbox{\it and likewise the curve C is also intersected by almost every "horizontal" line.\  }\end{eqnarray*}

Let $\tilde C\in \Bbb C^2$ be an analytic curve with the same property. Then we can construct a John mapping $T$ in the same way as in the real case. And we can ask a similar question about {\it the periodicity problem for the complex John mapping} :
What curve from a given class of curves has the following property?
\be\label{CJPerCond}\quad   The \ comp\,lex\ John\ mapping\ \
T:\tilde C\to \tilde C\  \ has\ at\ least\ one\ periodic\ point.\ee

In accordance with a corresponding complex setting of the problem \re{SE},\re{NDPH}, we will use the following:
\begin{eqnarray}\label{ComlSet}& \mbox{ \it
To find two meromorphic functions }\ f(x),\,g(y)\ \mbox{\it of one
complex variable such} \nonumber \\ & \mbox{\it that the condition }\
f(x)+g(y)=0\ \mbox{\it is satisfied as soon as} (x,y)\in \tilde C\subset\Bbb
C^2.\ \ \end{eqnarray}

Note that for the case of transitive action of the mapping $T$, the proof of the proposition 1 can easily be extended to this complex case:
\begin{pr}\label{uniqueCom}
For a biquadratic complex curve \re{2_2F}, if the complex John
mapping $T$ transitively acts on $\tilde C$, then the complex problem
\re{ComlSet} has only a trivial solution.
\end{pr}

Below we will consider biquadratic curves $C$ satisfying the property \re{CoDom} and we will give an explicit criterion, distinguishing the cases of periodicity and transitivity. In the former case, we will construct an explicit nontrivial solution of the problem \re{SE},\re{NDPH} in sense of
the settings \ref{UHDP} or \ref{MHDP} and as an intermediate setting \re{ComlSet}. In the latter case, i.e. in the case of transitivity, every solution will be proven to be a zero solution.

\subsection{Boundary value problems and the moment problem}\label{DirNeuMo}

In this subsection we intend to indicate equivalence of properties of some boundary value problems, in particular, of the Dirichlet and Neumann problems, and present a moment problem which is responsible for these properties. Thus we will obtain some problems that are equivalent to the problem \re{SE},\re{NDPH} in the setting \re{UHDP}. For details and generalizations, see \cite{BurMomProblem} or \cite{BurBook}. This equivalence is based on the connection condition of solution traces for the Cauchy problem.

Let us consider a hyperbolic equation in an arbitrary bounded
domain $\Omega\subset\Bbb R^2$ with a smooth boundary \be L u = a
u''_{x_1 x_1} + b u''_{x_1 x_2} + c u''_{x_2 x_2} = 0 \lab{GE}\ee
in the Sobolev space $H^m (\Omega), m \ge 3$.

In addition, let us introduce a conormal vector ?¤ and a derivative with respect to the conormal by means of an analog of the Green formula for the Laplace operator
$$\int_\Omega (L u \cdot \overline v - u \cdot L
\overline v)\; d x = \int_{\partial \Omega} (u^\prime_{\nu_*}
\overline v - u \overline v^\prime_* )\; d s. 
$$

One can consider that $\frac {\partial}{\partial \nu_*} = l (\nu)
\frac {\partial}{\partial \nu} - \frac 1{2 k} [ l (\nu (s))
]^\prime_s \cdot \frac {\partial}{\partial s},$ $l (\xi)
=a\xi_1^2+b\xi_1\xi_2+c\xi_2^2\ $ is the symbol of the operator $L
$, $ \nu $ is a unit vector of the normal, $s$ is a natural parameter
on $\partial \Omega$, $k=\pm|\nu'_s|$ is the curvature or,
more precisely, $\nu'_s=k\tau$ where $\tau=(-\nu_2,\nu_1)$ is the tangent
vector.

And let us consider an over-determined boundary value problem for the equation \re{GE} \be u^\prime_s \vert_{\partial \Omega} =
\gamma, \,\, u^\prime_{\nu_* }|_{\partial \Omega} = \kappa.
\lab{C*P}\ee
which is can also be written as the Cauchy problem $u
\vert_{\partial \Omega} = \gamma_0, \,\, u^\prime_{\nu
}|_{\partial \Omega} = \kappa_0.$

The following question can be formulated here, which is: what is a connection between the functions $\gamma$ and $\kappa$ if they are generated by solution $u$ of the equation \re{GE}? In order to answer the question, we will need a certain construction, which we are going to describe now.

The equation \re{GE} can be rewritten also as
\be (\nabla \cdot a^1) (\nabla \cdot a^2) u = 0\lab{ENabla}\ee
where $a^j = (a_1^j, a_2^j), j = 1, 2 $ are some unit real vectors.
Let us introduce vectors $ \tilde a^1 = (-a_2^1, a_1^1), \, \, \tilde
a^2 = (-a_2^2, a_1^2)$, which are some direction vectors of a set of characteristic
directions $ \, \, \Lambda = \Lambda^1 \cup \Lambda^2 $, $
\Lambda^j = \left \{\lambda \tilde a^j \vert \lambda \in {\Bbb R}
\right \} $, $ \, \, j = 1, 2, $ $ \langle \tilde a^j, a^j \rangle
= 0 $. Let us introduce also an angle $\varphi_0=\varphi_1-\varphi_2$, in which $\varphi_j$ is
any solution of equation $\tan \varphi_j=\lambda_j$, i.e.
$\varphi_j$ is an inclination angle of a vector of the
characteristic direction corresponding to the root $\bf \lambda_j
$, $\varphi_0$ is angle between characteristics, and let
$\Delta=\sin \varphi_0=\det\Vert a^1\; a^2\Vert,$ here $a_j$ are
columns. Here and below the vector of the characteristic direction
must be understood as a vector $ \nu\in
{\Bbb C \,} ^ 2 $ which is a null of the 
symbol: $ l(\nu) =0. $ 
The traces $\gamma$ and $\kappa$ of solution $u$ are linked by
the following relations.

\begin{st}
If a function $u \in H^m (\Omega), $ $m > 3 $ is a solution of the
problem \re{C*P} for the equation \re{GE} then the functions $
\gamma \in H ^ {m - 3 / 2} (\partial \Omega) $, $ \, \, \kappa \in
H ^ {m - 3 / 2} (\partial \Omega) $ from \re{C*P} satisfy the
conditions \be \forall Q \in {\Bbb R\,} [z], \, \int_{\partial
\Omega} \left [\kappa(s) + \frac {\Delta}2 \gamma (s)\right ] \,
Q(x(s) \cdot \tilde a^1)\, d s = 0, \lab{uu_*1} \ee \be \forall Q
\in {\Bbb R\,} [z], \int_{\partial \Omega} \left [ \kappa(s) -
\frac {\Delta}2 \gamma(s) \right ] \, Q(x(s) \cdot \tilde a^2)\, d
s = 0, \lab{uu_*2} \ee where $x(s)$ is a moving point on $\partial
\Omega $. \end{st}

An inverse statement holds also:

\begin{st}\label{statem2}
If functions $ \gamma \in H ^ {m - 3 / 2} (\partial \Omega) $, $
\, \, \kappa \in H ^ {m - 3 / 2} (\partial \Omega),  m>3 $ satisfy
the conditions \re{uu_*1},\re{uu_*2} then there exists a solution
$u \in H^{m-1-\epsilon} (\Omega)$ of the problem \re{C*P} for the
equation \re{GE} with each $\epsilon
>0.$ For functions
$\psi=u|_{\partial\Omega}, $ $\chi=u'_\nu|_{\partial\Omega}$ 
we have $l (\nu) \chi = \kappa + [l (\nu)] ^ \prime_\tau/2k\,\gamma $ and
also $ \, \, \psi = \int \gamma (s) d s + \mbox {const}  $  (the
Luzin's trigonometrical integral). In addition, function $u $ is restored univalently up to 
an additive constant. The mapping: $$H^{m} (\partial \Omega)
\times H ^{m} (\partial \Omega)/\{\mbox{const}\} \ni \{(\gamma,
\kappa) \mbox{ with \re{uu_*1},\re{uu_*2}}\} \to u \in
H^{m-1-\epsilon} (\Omega) $$ is continuous (for all $\epsilon>0$).
\end{st}

\begin{cor}
For each solution $u \in H^m (\Omega), m> 3 $ of the equations
\re{GE} the following Zhukovsky's equality \be \int_{\partial
\Omega} \kappa d s = 0 \lab {ZhC} \ee holds. \end{cor}

\begin{proof} \ It follows from the condition \re{uu_*1} for $Q \equiv 1 $
because $ \int _ {\partial \Omega} \gamma (s) d s = 0 $.
\end{proof}

\medskip Consider the following moment problem:
$$
\int _ {\partial \Omega} \alpha (s) (x (s) \cdot \tilde a^j) ^N d
s = \mu^j_N; \, \, j = 1, 2; \, \, N \in {\Bbb Z} _ +,
$$ where on the two given vectors $ \tilde a^j \in {\Bbb R \,} ^ 2 $
and on the two sequences of numbers $ \mu_N^j $
the function $ \alpha $ can be obtained.

Obviously, for the case when $\partial
\Omega$ is the unit circle and vectors $\tilde a^j,\,j=1,2$ are
equal $\tilde a^1=(1,i);\tilde a^2=(1,-i)$
this moment problem turn on well-known trigonometric moment
problem because then $(x (s) \cdot \tilde a^j) ^N =\exp(\pm iN)$.
Another way to the same is to write the Chebyshev polynomial $T_N$
insteed of the power.

Among numerous problems, connected to the above moment problem, we will consider the problem of indeterminacy (uniqueness), which can be formulated as follows:
for what curve $\partial \Omega$ and vectors $\tilde a^j,\,j=1,2$ can there
exist function $ \alpha $ such that \be \forall N \in {\Bbb
Z}_+, \, j=1,2,\, \int _ {\partial \Omega} \alpha (s) (x(s) \cdot
\tilde a^j) ^N d s = 0. \lab {HMP} \ee

The following fact is valid here:

\begin{st}\label{DirNeuMom}
Let $m \ge k\ge 3 $ and let us consider three sets of statements:

$1_m $) The homogeneous moment problem \re{HMP} has a nontrivial
solution $ \alpha \in H ^ {m - 3 / 2} (\partial \Omega) $.

$2_k $) The Dirichlet problem $u \vert _ {\partial \Omega} = 0 $
for the equation \re{GE} has a nontrivial solution $u \in H^k
(\Omega) $.

$3_k $) The Neumann problem $u ^\prime _ {\nu_ *} \vert _
{\partial \Omega} = 0 $ for the equation \re{GE} has a non-constant
solution $u \in H^k (\Omega) $.

Then $1_m) \Rightarrow 2 _ {m - q} $); $ \, \, 1_m) \Rightarrow 3
_ {m - q}); $ $ \, \, 2_m) \Rightarrow 1_m); $ $ \, \, 3_m)
\Rightarrow 1_m) $ with $q=1+0$ (By definition, for bounded domain
$H^{k+0} (\Omega) =\bigcup_{\epsilon > 0} H^{k+\epsilon}
(\Omega)$).
\end{st}

\begin{proof}
1) $ \Rightarrow $ 2). Using the pair $ \gamma = 0 $, $ \, \, \kappa =
2 \alpha / \Delta $, with the help of statement 1, we consrtuct the
solution $u \in H^{m-q} (\partial \Omega) $.

2) $\Rightarrow $ 1). We put $ \alpha = \kappa $ and apply the
statement \ref{statem2}.

The implications 1) $ \Rightarrow $ 3) and 3) $ \Rightarrow $ 1)
are similar. \end{proof}

Note that instead of considering the Neumann problem in statement 3, we can write down the boundary condition of the form
\be u_{\nu_*}=\lambda
u_{\gamma}\lab{GBC}\ee with an arbitrary constant $\lambda$.

Note also that the case, in which the domain $\Omega$ is an ellipse, was studied in the works of A. Huber \cite{Hub}, R. Alexandrjan
\cite{Alex2}, V.I. Arnold \cite{Arnold}, and, likewise, by one of the authors of this current paper in the works \cite{BurMomProblem} or \cite{BurBook}. An answer to the question regarding properties of such a Dirichlet problem as described by \re{SE},\re{NDPH} can be formulated through the following. First, let us reduce our problem considered inside an ellipse, by means of a linear transform, to the problem \re{NDPH} in a unit disk for the equation \re{ENabla}. Find slope angles $\varphi_1,\varphi_2$ of characteristics and an angle
$\varphi_0=\varphi_1-\varphi_2$ between the characteristics.

\begin{st} (\cite{BurMomProblem}, see also section \ref{breakdown})
The problem \re{SE},\re{NDPH}
has a nontrivial solution in the unit disk in a space
$H^k(\Omega),k\ge 2$ if and only if  \be\varphi_0/\pi\in \Bbb
Q.\lab{NC}\ee If the condition \re{NC} is satisfied, then there exists a denumerable set of linear independent polynomial solutions of the problem  \re{SE},\re{NDPH}.
 \end{st}

Thus, existence of a nontrivial solution of the Dirichlet problem \re{SE},\re{NDP} in a general bounded domain with a smooth boundary is equivalent to the existence of a non-constant solution of the boundary value problem \re{SE},\re{GBC}, in particular, of the Neumann problem, and it is equivalent to existence of a nontrivial solution of the moment problem \re{HMP}. Below we will provide a criterion of nontrivial solvability of each of these problems with the curve \re{2_2F}, in view of the condition \re{NC}.

---------------------------------

\subsection{Equation-domain duality and one more equivalent problem} 

Let $\Omega\subset{\Bbb R}^n$ be a bounded semi-algebraic domain
given via the inequality $\Omega=\{x\in{{\Bbb
R}}^n:P(x)>0\}$ with a real polynomial $P$. Equation $P(x)=0$
gives the boundary $\partial \Omega.$ We assume that the boundary of
the domain $\Omega$ is nondegenerate: $\ \mid\nabla P\mid\ne 0$ onа
$\partial\Omega.$ Consider the Dirichlet boundary value problem
for the equation \re{GE} of order $2$ with constant complex
coefficients: \be Lu=L(D_x)u(x)=0,\quad u\mid_{\partial\Omega}=0,
\lab {BBE}\ee where $D_x=-i\partial/\partial x$. We understand the
equation-domain duality as a correspondence 
the problem \re {BBE} and the equation:
\be P(-D_\xi)
\{L(\xi)\,w(\xi)\}=0,\lab {DBE}\ee given in the following
statement:

\begin{st}\label{Duality}
For each nontrivial solution of the problem \re {BBE} from $
C^2(\overline\Omega)$, there exists a unique nontrivial, analytic in
${\bf C}^n$, solution of the equation \re {DBE} from a cercain class $Z$
of entire functions and, conversely: for each nonzero solution
$w\in Z$ of the equation \re {DBE}, there exists a nonzero solution
$ u\in C^2(\overline \Omega)$ of the problem \re {BBE}. The class
$Z$ is determined here as space of Fourier images of functions of
the form $\theta_{\Omega}v;$  where $v\in C^{2}({\Bbb R}^n)$ and
$\theta_{\Omega}$ is the characteristic function of the domain
$\Omega.$ \end{st}

For clarity, we will provide here { \it a sketch of the proof.} Let us assume
that the problem \lab {BE} has a nontrivial solution $u$ in $C^{2}(\overline \Omega)$, and also let
$\tilde u\in C^{2}(\Bbb R^n)$ be any smooth continuation of $u$ on the $\Bbb R^n$; then we multiply $\tilde u$
by the characteristic function $\theta_{\Omega}$ of domain $\Omega$ ($\theta_{\Omega}=1$ in domain
$\Omega$ and ($\theta_{\Omega}=0$ outside of $\Omega$) and then apply the operator $L$ to the product $\theta_{\Omega}\tilde u$.

Differentiating the product
by the Leibniz rule, we obtain \be L(\theta_{\Omega}\tilde
u)=\theta_{\Omega}L(\tilde u)+L_1(u,\nabla
u)\delta_{\partial\Omega}+L_2(u)(\delta_{\partial\Omega})^\prime_{\nu}\lab{bbb}\ee
where $\delta_{\partial\Omega}$ is some measure supported on
$\partial\Omega$:
$<\delta_{\partial\Omega},\phi>=\int_{\partial\Omega}\overline\phi(x)
ds_x,$ $L_2(u)=L(\nu)u,$ $ \nu$ is the external normal as before.

The first term in \re{bbb} is equal to zero by means of the equation.
Taking into account the boundary condition $u_{\partial\Omega}=0$,
one will have the last term being transformed to the following term
(in view of the second term):$L_1(u,\nabla u)\,\delta_{\partial\Omega}$
as, for example, in the case of one variable $x\delta'(x)=-\delta.$
Then, multiplying the obtained equality by $P(x)$, the right-hand
side vanishes because, for example, for the case of one variable $x\delta(x)=0.$
We obtain an equation of the form $P(x)L(D)(\theta_{\Omega}\tilde u)=0.$
Application of the Fourier transform $\mathcal F$ leads to the equation
\re {DBE} where $w=\mathcal F(\theta_{\Omega}\tilde u).$
Hence, necessity is proven. Sufficiency will be obtained by means of
conversion of this proof. For a full proof for a general case,
see the works \cite{BurBook} or \cite{BurAlgB}.

The term {\it'equation-domain duality'}\ here has a meaning of equivalence
of problems \re {BBE} and \re {DBE}, which reads here as follows:

\be \qquad L(D_x)u=0,\quad
u\mid_{P(x)=0}=0,\qquad \lab {BEE}\ee \be P(-D_\xi)v=0,\quad
v\mid_{L(\xi)=0}=0. \lab {DBEE}\ee

---------------------------------------------------

\noindent
Now, from statement \ref{Duality}, it follows:

\begin{pr}
If there exists a nontrivial solution of the problem \re{BEE} with
the string equation \re{SE}, $L(\xi)=-\xi^1\xi^2$ where
$\xi=(\xi^1,\xi^2)$ is a covector, in a plane bounded domain with
the biquadratic curve \re{2_2F} as a boundary $\partial\Omega$,
then there exists also a nontrivial solution $v\in Z$ of the following
problem \be\sum\limits_{i,k=0}^2
a_{ik}\frac{\partial^{i+j}v}{\partial \xi^i\partial
\eta^k}=0,\quad v|_{\xi=0}=0,\ v|_{\eta=0}=0.\lab {DBEEE}\ee
The opposite holds as well.
\end{pr}

The last boundary value problem for the equation of the fourth order
has only two boundary conditions instead of four as it would be,
for example, for the problem of the Goursat (or Darboux) type.
Therefore it is not surprising that a nontrivial solution of the
problem \re {DBEEE} exists here. But, as will be shown below, for almost
every curve \re{2_2F}, the problem \re{SE}, \re{NDP} has only a trivial
solution; therefore, almost every problem \re {DBEEE} has only
a trivial solution as well. But this statement can seem
to be astonishing, namely, due to insufficiency of data,
in spite of a stipulation that the solution $v$ belongs to the space $Z$.

\vskip 5pt
We finished stating propositions on boundary value problems
for general domains and now we have to wait till we obtain
explicit answers for domains with biquadratic boundaries.

\section{Generic biquadratic curve}

\subsection{Parameterizations of generic biquadratic curve} 


The complex curve \re{2_2F} is remarkable, as it is the most general
algebraic curve having the following property: almost every vertical
or horizontal line, which can intersect $C$, will intersect $C$ in 2 points.

Let \re{2_2F} be a generic nondegenerate real biquadratic curve.
Assume that the parameters $a_{ik}$ are chosen such that the real
curve $C$ bounding the domain satisfies the condition \re{CoDom}
of the subsection \ref{JohnC}.

We will begin our study of the problem (1.1),(1.3) by an observation
that the curve (1.4) is an elliptic curve allowing for uniformization
in terms of elliptic functions.

Indeed, the equality \re{2_2F} can be rewritten in one of the two forms
\be
A_2(x)y^2 + A_1(x) y +A_0(x) =0 \lab{A_form} \ee or  \be B_2(y)x^2
+ B_1(y) x +B_0(y) =0, \lab{B_form} \ee

where $A_i(x)$ and $B_i(y)$ are polynomials of order two or lower.
Multiplying the equation \re{A_form} by $A_2$ and the equation
\re{B_form} by $B_2$, these expressions reduce to the forms
\be Y^2-D_1(x) =0 \quad \mbox{or}
\quad X^2-D_2(y) =0, \lab{full_sq} \ee where $Y=2A_2(x)y+A_1(x),
\; X=2B_2(y)x + B_1(y)$; $\ (x,y)\to(x,Y)$, $(x,y)\to(X,y)$ are
birational transformations and $D_1(x), D_2(y)$ are
discriminants of quadratic equations \re{A_form} and \re{B_form}:
$$
D_1(x) = A_1^2(x) - 4 A_0(x) A_2(x), \; D_2(y) = B_1^2(y) - 4
B_0(y) B_2(y).
$$

In a general situation, the discriminants $D_{1,2}$ are polynomials
of order 4 or 3. From the work \cite{WW}, it follows that every curve of the kind
\be Y^2 = \pi_4(x) \lab{gen_e4} \ee
with a generic fourth order polynomial $\pi_4(x)$ can be transformed to the following canonical form
\be Y^2 = 4 x^3 - g_2 x - g_3 \lab{canonical_3} \ee
which allows a standard parameterization
\be x= \wp(t) , \; Y= \wp'(t)
\lab{par_wp} \ee
through the elliptic Weierstrass function $\wp(t)$ with primitive periods
$2\omega_1,\ 2\omega_2$. The parameters  $g_2, g_3$ are so-called (relative)
invariants of the polynomial $D_{1}$. They are real because $\pi_4$ is real.

As the transformation $(x,Y)\to(X,y)$ is birational, primitive periods
of the curves \re{full_sq} coincide. The invariants $g_2,g_3$ can be found
through periods. Hence we obtain the following important statement
(mentioned by Halphen \cite{Halphen}):

\begin{st}\label{invar} 
Invariants $g_2,g_3$ of polynomials $D_1(x)$ and $D_2(x)$ are
identical.
\end{st}

Thus, both curves \re{full_sq} are elliptic curves (see \cite{WW}) with appropriate primitive periods $2\omega_1, 2\omega_2$, the curve \re{2_2F} is homeomorphic to the torus: $C\approx \Bbb
C/(2\omega_1\Bbb Z\oplus 2\omega_2\Bbb Z)$ and there exist standard structures of a commutative group and an Abelian variety \cite{GHP}. We deal with elliptic functions of the second order. Recalling \cite{WW}, it can be noticed that properties of a general elliptic function can be characterized by the number of poles (or, equivalently, zeroes) in the fundamental parallelogram of periods. This number is called an order of the elliptic function (the order is taken with the account of multiplicity of poles). From the Liouville theorem it follows that the simplest possible order is two \cite{WW}. For example, the Weierstrass function $\wp(t)$ is of order 2 because it has only one pole and that pole is of multiplicity 2 at the point $t=0$ of the fundamental parallelogram.

The general elliptic function of the second order $\Phi(t)$ has
two arbitrary poles $p_1,p_2$ and two arbitrary zeroes $\zeta_1,
\zeta_2$ in the parallelogram of periods. The only condition is
$p_1+p_2-\zeta_1-\zeta_2=\Omega$, where $\Omega =2m_1
\omega_1 + 2 m_2 \omega_2$ is an arbitrary period (\cite{WW}). It can be
easily shown 
that a generic elliptic function of the
second order with given periods $2 \omega_1, 2 \omega_2$ can be
presented as \be \Phi(t) = \frac{\alpha \wp(t-t_0)+\beta}{\gamma
\wp(t-t_0) + \delta} \lab{Phi_wp} \ee Thus, $\Phi(t)$ depends on four
independent parameters, for example, $\alpha, \beta, \delta, t_0$.

There may be another, sometimes more convenient, representation of the
function $\Phi(t)$: \be \Phi(t) = \kappa \: \frac{\sigma(t-e_1)
\sigma(t-e_2)}{\sigma(t-d_1)\sigma(t-d_2)}, \lab{Phi_sig} \ee
where $\kappa$ is a constant and parameters $e_1,e_2,d_1,d_2$ are
related to each other via \be e_1+e_2=d_1+d_2. \lab{p_z} \ee

The form
\re{Phi_sig} is obtained from a standard representation of an
arbitrary elliptic function expressed via the Weierstrass
sigma-functions \cite{WW}. In this case the points $e_1,e_2$ and
$d_1,d_2$ coincide with zeros and poles of the function $\Phi(t)$
and the condition \re{p_z} is equivalent to a balance condition
between poles and zeros of the generic elliptic function.

Note that apart from $\wp(t)$ there are another special cases of
functions of second order, such that the Jacobi elliptic functions
$sn(t;k), cn(t;k),dn(t;k)$ \cite{WW}, which we will use below.

What is uniformization for the biquadratic curve? The answer can
be given through the following:
\begin{tm}\label{GBiqC}
The generic complex biquadratic curve \re{2_2F} can be
parameterized by a pair of elliptic functions of the second order
having identical periods: \be x(t) = \Phi_1(t), \quad y(t) =
\Phi_2(t) \lab{xy_Phi} \ee Conversely, any two elliptic functions
$x=\Phi_1(t), y=\Phi_2(t)$ of the second order having identical
periods satisfy an equation \re{2_2F}.
\end{tm}

This theorem was proven essentially, for example, in Halphen's
famous monograph \cite{Halphen} on elliptic functions.
We will give a proof based on some of Halphen's major ideas here.

\vskip10pt
{\it Proof.}  From \re{2_2F}, the Euler
differential equation (\cite{Ince}) follows in the form \be \frac{dx}{\sqrt{D_1(x)}} =
\pm \frac{dy}{\sqrt{D_2(y)}} \lab{el_DD} \ee because
$dy/dx=-F_x/F_y=-X/Y$ (see \re{full_sq}).

As we already saw, the polynomials $D_1(x)$ and $D_2(y)$ have identical
invariants $g_2,g_3$. Hence, they both can be
reduced to the same canonical Weierstrass form (see, e.g.
\cite{Akhiezer}, s.34) by means of a pair of the M\"obius
transforms \be\tilde x ={\mu_1 x + \nu_1\over\xi_1 x + \eta_1},
\quad \tilde y = {\mu_2 y + \nu_2\over\xi_2 y + \eta_2},\quad
\mu_i\eta_i-\nu_i\xi_i=1,\ i=1,2 \lab{MT}\ee with some complex
parameters $\mu_1, \dots \eta_2$. Hence the equation \re{el_DD}
becomes  \be \frac{d \tilde x}{\sqrt{4 \tilde x^3 -g_2 \tilde x -
g_3}} = \frac{d \tilde y}{\sqrt{4 \tilde y^3 -g_2 \tilde y -
g_3}}\, \lab{el_DDW} \ee because $d\tilde y/d\tilde x=-\tilde
X/\tilde Y$ as above. But the equation \re{el_DDW} implies that for
appropriate periods $2\omega_1,2\omega_1$\be \tilde x=\wp(u),
\quad \tilde y=\wp(u+u_0), \lab{wp_par} \ee where
$\wp(u)=\wp(u;\omega_1,\omega_2)$, $u$ is an uniformization
parameter and $u_0$ is a complex constant.

Now we can return to initial variables $x,y$ via using inverse
M\"obius transforms to to finally arrive at \be x=\frac{\alpha_1 \wp(u) +
\beta_1}{\gamma_1 \wp(u) + \delta_1}; \quad y=\frac{\alpha_2
\wp(u+u_0) + \beta_2}{\gamma_2 \wp(u+u_0) + \delta_2},
\lab{wpar_xy} \ee
where constants $\alpha_ i, \dots,
\delta_i$ satisfy the equation: $\quad \alpha_i\delta_i-\beta_i\gamma_i=1,\ i=1,2$.
Then the formula \re{Phi_wp} is applied which means that we,
in fact, obtained a pair of elliptic functions of the second order.

The inverse statement of the theorem follows from a general
theorem stating that any two elliptic functions $x(t)$ and $y(t)$ with the
same periods satisfy an algebraic equation $F(x(t),y(t))=0$. The
degrees of the polynomial $F(x,y)$ with respect to variables $x$
and $y$ are determined from orders of corresponding elliptic
 functions. If both functions are of order two, then the polynomial
 $F(x,y)$ has, at most, the second degree with respect to each variable
 (this statement is included in \cite{WW}, namely into a statement of one of the problems offered to a reader to solve).$\Box$

Thus, we proved the theorem \ref{GBiqC}. Moreover, simultaneously the following has just been proven

\vskip 5pt
\begin{pr} \label{LFT} There exist complex linear-fractional transforms \re{MT} such
that transformed generic curve, which has the same form \re{2_2F}, can
be parameterized only through the Weierstrass function $\wp(u)$ as in
\re{wp_par}.

\end{pr}

Taking into a consideration the expressions \re{xy_Phi}, our
parameterization can be expressed through the Weierstrass sigma-function as: \be
x(t) = \kappa_1 \: \frac{\sigma(t-e_1)
\sigma(t-e_2)}{\sigma(t-d_1)\sigma(t-d_2)}, \quad y(t) = \kappa_2
\: \frac{\sigma(t-e_3) \sigma(t-e_4)}{\sigma(t-d_3)\sigma(t-d_4)}
\lab{xy_gen} \ee with two restrictions $e_1+e_2=d_1+d_2$ and
$e_3+e_4=d_3+d_4$.

If now the shift $t \to t-t_0$ is performed of the uniformization
parameter $t$, then $t_0$ can be chosen such that, say,  $t \to t-t_0$. This means

\begin {pr} In the expressions  \re{xy_Phi} the function $\Phi_1(t)$ can be chosen to be even:
$\Phi_1(-t) = \Phi_1(t)$ through choosing the shift  $t \to t-t_0$.
\end{pr}

From our considerations, an important corollary (also mentioned by Halphen [31])
follows in the form:

\begin{pr} Consider the differential equation \re{el_DD}.
Let $D_1(x)$ and $D_2(y)$ be polynomials of degree 4 or 3 with identical
invariants $g_2,g_3$. And also let $(x(t), y(t))$ be a solution of
this equation (parameterized, for example, by an initial condition). Then
$x(t)$ and $y(t)$ satisfy a biquadratic equation of the form
\re{2_2F}.\end{pr}

\subsection{Biquadratic foliation and singular points}\label{BiFol}

There is an interesting mechanical interpretation of these last results.
Assume that we have a dynamical Hamiltonian system for two
canonical variables $x(t),y(t)$ satisfying a system of equations
\be \dot x = \frac{\partial H(x,y)}{\partial y}, \quad \dot y =
-\frac{\partial H(x,y)}{\partial x}, \lab{Ham_xy} \ee where
$H(x,y)$ is a Hamilton function of the system. Obviously $H(x,y)$
is an integral of the system \re{Ham_xy}, i.e. $
(H(x,y))^\cdot=0$. Let us choose the Hamiltonian as the biquadratic
function \re{2_2F}: $H(x,y)=F(x,y)$. Then $F(x,y)=c$ with some
constant $c$ that depends on initial conditions for $x$ and $y$. This
constant can be incorporated into the coefficient $\ a_{00}$, so we
can write down $\tilde F(x,y)=0$, where $\tilde F(x,y)=F(x,y)-c$
is again a biquadratic curve (note that for $\tilde F(x,y)$ the
coefficients $A_2(x), A_1(x), B_2(y), B_1(x)$ remain intact,
whereas the coefficients $\tilde A_0(y)$ and $\tilde B_0(y)$
change from their initial values by a constant). Then

$$\frac{\partial H(x,y)}{\partial y}=2 A_2(x) y + A_1(x) = \pm
\sqrt{\tilde D_1(x)},$$ where $\tilde D_1(x) = A_1^2(x) - 4 A_2(x)
\tilde A_0(x)$ (as $y$ can be excluded as a root of the quadratic
equation $A_2(x) y^2+A_1(x) y +\tilde A_0(x)=0$). Similarly
analogously $$\frac{\partial H(x,y)}{\partial x}=2B_2(y) x +
B_1(y) = \pm \sqrt{\tilde D_2(y)}.$$ Thus, we can observe that for any
fixed Hamiltonian level $H=c$, the variables $x(t)$ avd $y(t)$ satisfy
the Euler equation
\re{el_DD} in which polynomials $\tilde D_1(x), \: \tilde D_2(y)$
have identical invariants $\tilde g_2$ and $ \tilde g_3$. Note, that in
this case the invariants $\tilde g_2, \tilde g_3$ (and hence the
periods $2\omega_1, 2\omega_2$) will depend on the value of the
integral $c$. Thus, we obtain a whole one-parameter family of
biquadratic curves $F(x,y)=c$ and corresponding elliptic functions
$x(t), y(t)$ of the second order, which are trajectories of this
dynamical system.
\vskip 5pt

Now, let us rewrite discriminants $\tilde D_1(x), \tilde D_2(y)$ in factorized
forms
$$
\tilde D_1(x) = q_1\prod_{i=1}^4(x- x_i), \quad \tilde D_2(y) =
q_2\prod_{i=1}^4(y- y_i),
$$
where $q_1,q_2$ are leading coefficients of the discriminants and
$x_i, y_i, i=1,2,3,4$ are their roots (for simplicity we assume
that both discriminants have the degree 4). In general, roots $ x_i,
y_i$ will depend on the parameter $c$. What is mechanical meaning
of points $x_i, y_i$? From equations of motion it can be seen that
$$
\dot x = \pm \sqrt{\tilde D_1(x)}, \quad \dot y = \pm \sqrt{\tilde
 D_2(y)}
$$
Thus, $ x_i$ and $ y_i$ are {\it stable points}: ${\dot {\tilde
x}_i} = \dot {\tilde y}_i =0$. We require that the points
$(x_i,y_k)$ belong to our biquadratic curve $\tilde F(x,y)=0$.
This leads to a situation, in which points $(x_i,y_k)$ satisfy the
following conditions \be \tilde F(x,y)=0, \quad
\partial_x \tilde F(x,y)=0, \quad \partial_y \tilde F(x,y) =0
\lab{stab} \ee which, in turn, means that the points $(x_i,y_k)$ are
singular points of the biquadratic curve $\tilde F(x,y)=0$. As we
saw, in a generic situation, this curve is elliptic and, hence, has a
genus 1 (genus $C=C_{n-1}^2-d$, where $n$ is a degree and $d$ is the number
of double points of $C$). The order of this curve is 4. Assume
that the curve is irreducible. Then
such curve cannot have more than three singular
points in a complex projective plane. The latter is defined by the
coordinates $(s_0,s_1,s_2)$ such that $x=s_1/s_0, y=s_2/s_0$. In
these coordinates, the equation for our curve becomes \ba &&a_{22} s_1^2
s_2^2 + s_0 s_1 s_2 (a_{12}s_1 + a_{21} s_2) +
s_0^2(a_{20}s_1^2 + a_{02}s_2^2 + a_{11} s_1 s_2) + \nonumber \\
&&s_0^3 (a_{10} s_1 + a_{01} s_2) + s_0^4 \tilde a_{00}
=0.\lab{hom_F} \ea

Elementary considerations show that two points $(0,1,0)$ and
$(0,0,1)$ of the projective plane are singular for any values of
parameters $a_{ik}$. Thus, only one singular point can exist in
each finite part of the plane.  In turn, this means that, at least,
two roots, say $x_1, x_2$ of the discriminant $\tilde D_1(x)$
should coincide: $ x_1= x_2$. The same condition holds for the
discriminant $\tilde D_2(y)$, i.e. $ y_1= y_2$ because invariants
$g_2,g_3$ of the discriminants $\tilde D_1(x)$ and $\tilde D_2(y)$
are identical. Then it can be elementary verified that the point $(
x_1, y_1)$ will be, indeed, a singular point of the biquadratic
curve. In principle, the second singular point can occur. But in
this case the genus of the curve will be -1 meaning that the curve
is reducible.

Thus, we can formulate the following
\begin{pr}\label{singP}
The irreducible biquadratic curve $\tilde F(x,y)=0$ has a
singular point in a finite part of the (complex) plane if and only
if the discriminant $\tilde D_1(x)$ (and, hence, $\tilde D_2(y)$ as
well) has a multiple zero $x_1$ (and, respectively, $y_1$). In this
case, the point $(x_1,y_1)$ is singular and unique.
\end{pr}

Note that this proposition can be reformulated in an equivalent form.
Indeed, the polynomial $D_1(x)$ has a multiple zero if and only if
its discriminant is zero. Thus, in order to find all singular
points of the curve, first,
the discriminant $D_1(x)$ of the biquadratic curve $F(x; y) = 0$
must be calculated and, then, the discriminant  $\Delta(D_1(x))$ of
the polynomial $D_1(x)$ is calculated (i.e. discriminant of discriminant).
Obviously, the following two discriminants coincid
$$
\Delta(D_1(x))=\Delta(D_2(y))
$$
because the invariants $g_2,g_3$ of the polynomials $D_1(x)$ and
$D_2(y)$ are identical (statement \ref{invar}). Condition
$\Delta(D_1(x))=0$ (or, equivalently, $\Delta(D_2(y))=0$) leads to
nonlinear equations for the coefficients $a_{ik}$. Under such
condition the biquadratic curve $F(x,y)=0$ has a genus $<1$, i.e.
it is either irreducible and has only one singular point, which is
in a finite part of the projective plane, or reducible:
$F(x,y)=\tau_1(x,y) \tau_2(x,y)$, where $\tau_1(x,y), \tau_2(x,y)$
are two polynomials, which are linear with respect to each argument $x,y$
(but $\tau_k(x,y)$, in general, are not linear functions, as they may
contain the terms like $xy$). Hence, we have obtained

\begin{pr}\label{genusP} The condition $
\Delta(D_1(x))\ (=\Delta(D_2(y)))\not=0 $ is necessary and
sufficient for the equality:\ \ genus $C=1$.
\end{pr}

\subsection{Case of a generic symmetric curve}

Above we considered a generic case when our biquadratic curve
$F(x,y)=0$ was non-symmetric.

Now we will assume that our curve is symmetric, i.e. $F(x,y)=F(y,x)$.
Equivalently, this means that the coefficient matrix $a_{ik}$ in
\re{2_2F} is symmetric: $a_{ik}=a_{ki}$. Obviously, in this case
both discriminants coincide $D_1(x)=D_2(x)$. From the Euler
differential equation \re{el_DD} it can be deduced that a parameterization
can be given by the following formula

\be x(u)= \Phi(u), \quad y(u) =
\Phi(u+u_0), \lab{par_sym} \ee where $u_0$ is a constant and
$\Phi(u)$ is an even function of the second order, i.e. $\Phi(t)$
can be presented in the form \be \Phi(u)= \frac{\alpha \wp(u) +
\beta }{\gamma \wp(u) + \delta} \lab{sym_Phi} \ee Thus, in the
symmetric case a parameterization is provided through an even
elliptic function of the second order. Conversely, any pair
$(x,y)=(\Phi(u), \Phi(u+u_0))$ with an arbitrary $u_0$ generates a
symmetric biquadratic curve via expressions \re{par_sym} because the
point $(y,x)=(\Phi(u+u_o), \Phi(u))=(\Phi(\tilde u), \Phi(\tilde
u+u_0)),\, \tilde u=-u-u_0$ also belongs to $C$. Note that the last
statement was attributed to Euler in the work
\cite{VesUs}.

Thus, we obtain the following
\begin{pr}\label{www}
The generic complex symmetric biquadratic curve \re{2_2F} can be
parameterized by means of an even elliptic function of the second order and
a shift similar to that given in \re{par_sym}.
Conversely, for every elliptic function $\Phi(u)$ (no matter
whether it is even or not) of
the second order and any shift $u_0$, the variables $x,y$ from
\re{par_sym} satisfy the equation \re{2_2F} with a symmetric
matrix $A$.
\end{pr}

This appears important

\begin{pr}
The generic complex biquadratic curve \re{2_2F} can be transformed
to a symmetric curve by means of linear-fractional complex changes
of variables \re{MT}. If the initial curve is real then it can be
transformed to a real symmetrical curve also (although a corresponding
transformation can be performed also for a case with complex coefficients).
\end{pr}

\begin{proof} As we showed in the theorem \ref{GBiqC} (see proposition
\ref{LFT}), a generic non-symmet\-ric curve can be transformed to a
curve described by \re{wp_par} by linear-fractional complex
changes \re{wp_par} which means that the curve-image is symmetric in
virtue of the last proposition. In the real case we note that, as
is well-known, the invariants $g_2,g_3$ of real polynomial
$D_1(x)$ are real, so that the differential equation \re{el_DDW}
is real and there exists its real solution $(x(u),y(u))$ that can
be extended to a complex domain and, therefore, satysfies the
symmetric equation \re{2_2F}. The last equation has irreducible
polynomial which must be real as $a_{22}=1$ because it is always
possible to choose real parameters $u_1,...u_8$ such that the vectors $V_i$,
which are composed of components $x^k(u_i)y^l(u_i),\ 0\le k,l\le 2, 0\le
k+l<4$ will be linear independent, so that the coefficients of our
polynomial satysfy a linear system of eight linear equations with real
coefficients $(V_{i})_j$ and the real right-hand side parts
$-x^2(u_i)y^2(u_i)$.
\end{proof}


\noindent
Note, that in the work \cite{IatrouRob} there is another proof
of the latter fact (see statement \ref{StatReduce} below).

Let us now consider the question how one can restore the polynomial $F$ if
the discriminants $D_1,D_2$ are known. For a symmetric case,
T. Stieltjes \cite{Stieltjes} proposed a neat explicit formula for
the polynomials $F(x,y)$ by means of a solution of the
differential equation \re{el_DD}. Assume that
$$
D_1(x)=D_2(x)=b_0 x^4 + 4 b_1 x^3 + 6 b_2 x^2 + 4 b_3 x + b_4.
$$
Then $F(x,y)$ can presented via the following determinant: \be   F(x,y)= \left
|
\begin{array}{cccc} 0 & 1 & -(x+y)/2 &
xy\\ 1 & b_0 & b_1 & b_2-2 C \\ -(x+y)/2 & b_1 & b_2 + C  & b_3\\
xy & b_2-2 C & b_3 & b_4 \end{array} \right | \lab{det_F} \ee The
function $F(x,y)$ is defined up to an arbitrary non-zero number
factor. It is assumed that the curve is nondegenerate, i.e. its
genus is 1. It may be true if and only if the determinant
$$
\Delta = \left |
\begin{array}{ccc} b_0 & b_1 & b_2-2 C \\ b_1 & b_2 + C  & b_3\\
b_2-2 C & b_3 & b_4 \end{array} \right |
$$
is nonzero. $C$ is an (arbitrary) integration constant
(\cite{Stieltjes}).  The Stieltjes formula is useful in the
problem of reducing of the arbitrary symmetric biquadratic curve
to some standard forms, as will be seen below.

\subsection{Canonical forms of biquadratic curve}\label{CanForms}
The general curve \re{2_2F} contains eight
free parameters. It is natural to transform this curve to the
form containing the smallest possible numbers of parameters.

First, note that under some arbitrary projective transformations of the form \be x
\to \frac{\xi_1 x + \eta_1}{\mu_1 x + \nu_1}, \quad x \to
\frac{\xi_2 x + \eta_2}{\mu_2 x + \nu_2} \lab{proj_xy} \ee with
complex parameters one can obtain an equation similar to the equation \re{2_2F}, in which the parameters $a_{ik}$ will be some altered parameters. This idea has already exploited in the proof of the theorem 2. As every projective transformation
\re{proj_xy} contains three independent parameters, it is
possible to reduce the total number of independent parameters
$a_{ik}$ to $8-6=2$. As our curve is elliptic, these free
parameters will have only the invariants $g_2,g_3$ under linear-fractional
transformations, performed separately for each of the variables $x$ and $y$. More explicitly, if we consider the projectivisation
$\sum\limits_{i_1,i_2,j_1,j_2=1}^2 a_{i_1i_2j_1j_2}
x_1^{i_1}x_2^{i_2}y_1^{j_1}y_2^{j_2},\ x=x_1, y=y_1$ of the initial
curve \re{2_2F} and its projective transformations of variables
$x$ and $y$ separately, then $g_2,g_3$ will be the only invariants (having weights equal to 4 and 6, respectively), but for the transformations of the
group $SL(2,\Bbb C)$, these invariants will be absolute.

Above we reduce the general case to the symmetric case so that we
restricted ourselves to considering a symmetric curve
$F(x,y)=F(y,x)$. So now we would like to find some canonical forms of
the curve that may contain only two parameters.

There are two obvious canonical forms which can be obtained for
the curve $F(x,y)=0$. These two forms correspond to two canonical
forms of elliptic integrals in the Euler differential equation
\re{el_DD}.

\vskip5pt

(I) The first form can be obtained if one reduces polynomial
$D_1(x)=D_2(x)$ to the canonical cubic Weierstrass form:
$$
D_1(x) = 4 x^3 - g_2 x - g_3.
$$
Such form can be always obtained through using an appropriate projective
transformation (with possible complex coefficients). Then, from the
Stieltjes formula \re{det_F}, we obtain the following expression (see also
\cite{Halphen}) \be F(x,y)=(xy + (x+y)w + g_2/4)^2 -
(x+y+w)(4xyw-g_3)=0, \lab{WP_curve} \ee where $w=C$ is an
arbitrary parameter. There is a parameterization of this curve through the Weierstrass elliptic function $$x(u)=\wp(u), \; y(u)
= \wp(u+u_0)$$ and $w=\wp(u_0)$, where $u_0$ is an arbitrary
constant. The discriminant can then be easily calculated through
$$
D_1(x)=(4w^3-g_2 w-g_3)(4x^3-g_2 x-g_3)
$$
If the parameter $w$ is such that $4w^3-g_2 w-g_3=0$, then $D_1(x)
\equiv 0$ and in this case the curve $F(x,y)$ is reducible:
$F(x,y)=\rho^2(x,y)$, where $\rho(x,y)$ is a polynomial of order one
with respect to each variable $x,y$. If $4w^3-g_2 w-g_3 \ne 0$,
then the curve is irreducible. The singular points in a finite
part of the complex projective plane appear only if $g_2^3=27
g_3^2$. This condition means that the polynomial $4x^3-g_2 x-
g_3=4(x-e_1)(x-e_2)(x-e_3)$ has a double root, say $e_1=e_2$. In
this case, the following can be taken:  $g_2=3\tau^2, \; g_3=\tau^3$,
where $\tau$ is some parameter. Then it can be easily verified that the point $x=y=-\tau/2$ will be the only finite singular point of the curve
$F(x,y)=0$. If $g_2^3-27 g_3^2 \ne 0$, then there exsist no singular
points in the finite part of the projective plane and the curve
\re{WP_curve} is irreducible and has a genus of one.

The curve \re{WP_curve} contains three parameters $w,g_2,g_3$. Assume
that $4 w^3 -g_2 w -g_3 =W \ne 0$. In this case, the curve is
irreducible. By the linear transformation of arguments $x \to
\alpha x + w, \; y \to \alpha y + w$, where $\alpha^3 = 1/W$, the
terms $x^2y, xy^2$ and $x^2+y^2$ can be eliminated. Then, the curve reduces to the form
\be F(x,y)=x^2y^2 - x-y + Axy +B=0 \lab{WP_red} \ee which contains only two independent parameters $A,B$. The curve is elliptic and
non-singular (i.e. it has the genus of one) if condition \be
\Delta(D_1)=B(A^2-4B)^2 + A(36B - A^2) -27 \ne 0 \lab{WAB} \ee
holds.

Thus we proved the following

\begin{pr}
The generic complex biquadratic curve \re{2_2F} can be transformed
to the canonical form  \re{WP_red} by means of the linear-fractional
complex changes of variables \re{MT}. If the initial curve is real,
then it can be transformed to a real form also (although
corresponding transformation can be also performed for complex coefficients).
\end{pr}

As was seen above, any irreducible biquadratic curve can be transformed to fit this form. However, here we need general projective
transformations \re{proj_xy} with complex parameters. This will be
true even for a case, in which all parameters $a_{ik}$ of the
biquadratic curve are real.

\vskip 5pt

(II) To begin with, let us recall the well-known Legendre transformation.
Assume that an elliptic curve is already in the form $Y^2=\pi_4(x)$
(see \re{full_sq}), where $\pi_4(x)$ is an arbitrary polynomial of
the fourth degree with real coefficients. Using only a
linear-fractional transformation $x=\Gamma(\widetilde x)$ with
real coefficients, it is possible to reduce this curve to the form
(with a new $Y$) where $\pi_4(x)$ contains only even degrees with respect to $x$:
$\pi_4(x) = \alpha x^4 + \beta x^2 + \gamma$ with some real
coefficients $\alpha, \beta, \gamma$ (see, e.g. \cite{Bat}).

Now let us assume that our biquadratic curve is symmetric and contains only
real coefficients $a_{ik}$. Then it is always possible to transform this
curve to a form, that excludes odd degree terms, i.e.
$a_{21}=a_{12}=a_{10}=a_{01}=0$. Indeed, taking $x=y$ in the
expression \re{2_2F} for $F$, we obtain the polynomial $F(x,x)$ of
fourth order, and through applying the change
$x=\Gamma(\widetilde x)
$ to the polynomial in a generic case, the equation
reduces the equation $F(x,x)=0$ to the form $\alpha\widetilde x^4+\beta\widetilde x^2+\gamma=0$.
Therefore the change with the same $\Gamma$ \be
x=\Gamma(\widetilde x),\ y=\Gamma(\widetilde y)\lab{Gamma}\ee
applied to \re{2_2F} gives an equation with a desired polynomial
(observation from \cite{IatrouRob}).

We then obtain the so called Euler-Baxter biquadratic curve \cite{Bax},
\cite{Ves}: \be F(x,y) = x^2 y^2 + a(x^2+y^2) + 2 b xy + c=0,
\lab{E_Bax} \ee where $a,b$ and $c$ are real parameters if the
initial curve is real and, obvious, they will be complex if the
initial curve was complex. The curve \re{E_Bax} plays a crucial
role in derivations of the addition theorem for the elliptic function
$sn(t)$ \cite{Akhiezer}. The curve also appeared in Baxter's approach to
8-vertex model in statistical mechanics \cite{Bax}.

We first analyze possible finite singular points of the curve
\re{E_Bax}. The discriminant of the curve \re{E_Bax} is equal to
\be D_1(x)= -ax^4-(a^2-b^2+c)x^2-ac=-a(x^4-\tilde b x^2+c);\
\tilde b=\frac{b^2-a^2-c}{a},\label{D_1}\ee $D_1(x)=D_2(x)$ and and in virtue of the proposition \ref{singP} the curve \re{E_Bax} has a finite singular point under the condition
$$\Delta(D_1(x))=16a^2c\,((a^2-b^2+c)^2-4a^2c)^2=0.$$
If $c=0$, then there is an obvious singular point $x=y=0$. This
singular point will be isolated for all 
values of the parameters $a,b$, except for the case $a \pm b=0$. If
$a\pm b=0$, then the curve $F=0$ becomes reducible. In any case, if
$c=0$, then the genus of the curve \re{E_Bax} will be less then 1 (namely, it will equal 0 for the generic case and $-1$ for the exceptional case, which corresponds to a reducible curve). Thus, in this case the curve is not elliptic and can therefore be parameterized by means of rational functions. Here, the change $x \to
\kappa /x, \: y \to \kappa /y$ leads to the form $F=x^2+y^2 + axy
\pm1$. 

Then, using the scaling transformations $x \to \kappa x, \: y \to
\kappa y$ in the real case, the coefficients $c\neq 0$ can be reduced
to $\pm 1$, depending on a sign of this coefficient.

the Thus, in real symmetric case, the following simplest forms can be written:

(i) $F=x^2y^2 + a(x^2+y^2) + 2b xy +1$;

(ii) $F=x^2y^2 + a(x^2+y^2) + 2b xy - 1$.


\noindent 
In the form (ii), one assumes that $a>0$ by virtue of substitution $x\to
1/x,y\to 1/y$.

The canonical biquadratic curve \re{2_2F} with $F$ from (i) or
(ii) will be called by us as the Euler-Baxter curve (similarly to  \re{E_Bax}).

\begin{pr}\label{SymReduce}
The generic complex biquadratic curve \re{2_2F} can be transformed
to the canonical form (i) by means of the linear-fractional complex
changes of variables \re{MT}. If the initial curve is real, then it
can be transformed to the real form (i) or (ii) (although the
corresponding transformation may include complex coefficients). By
applying only the real linear-fractional transformation \re{MT},
any generic real symmetric biquadratic curve \re{2_2F} can be
reduced to one of the two forms: (i), (ii).
\end{pr}

Note that in the paper \cite{IatrouRob}, the following was proven

\begin{st}\label{StatReduce}
By applying only a real birational transformation, any generic real elliptic biquadratic curve \re{2_2F}
can be reduced to one of the three forms: (i), (ii)
or the form

(iii) $F=x^2y^2 + a(x^2-y^2) + 2b xy - 1$.
\end{st}


Note that the discriminants of the curve (iii) are equal to \be
D_1(x)= -a(x^4-\hat b x^2+1);\ D_2(y)= a(y^4+\hat b y^2+1);\ \hat
b=\frac{b^2+a^2+1}{a}.\label{D_1D_2}\ee

In the complex case, the coefficients $c\neq 0$ can be reduced to 1
and, hence, only the case (i) woll be left among the nondegenerate cases.
Following the work \cite{Bax}, one can find here a parameterization of the curve \re{E_Bax} through the elliptic Jacobi function $sn$:
\begin{pr}\label{Param} The curve \re{E_Bax} can be parameterized
through the expressions
\be x=\sqrt{k} \, sn(t;k), \quad y=\sqrt{k} \:
sn(t\pm\eta;k), \lab{Bax_par} \ee for any sigh $\pm$, where
parameters $k, \eta$ are determined through the following relations
\be k + k^{-1}
=(b^2-a^2-1)/a, \quad 1+ k\: a\: sn^2(\eta;k) =0 \lab{k_et}. \ee
\end{pr}

Note that one can perform the following substitution
$t= \tilde t+K$ to allow dealing with an even function.
Note that the mapping \re{Bax_par} is an analytic
diffeomorphism from a fundamental parallelogram, factored into a torus by the standard way, onto $\tilde C$.

Real parameterizations of the curves (i), (ii) and (iii) were given
in \cite{IatrouRob} and \cite{IatrouArx}.

Now let us recall the expressions \re{full_sq} and rewrite them as \be
F_x^2(x_1,y_1)=D_2(y_1)=B_1^2(y_1) - 4 B_0(y_1)
B_2(y_1)\lab{F_xD2}.\ee
For deducing the $y$-extreme point $(x_1,y_1)$, we will consider the following conditions: $0=y'(x_1)=F_x/F_y(x_1,y_1)$;
$0=F_x(x_1,y_1)=2B_2(y_1)x_1+B_1(y_1)$. If $B_2(y_1)=0$, then from \re{B_form} it follows that $x_1=-B_0(y_1)/B_1(y_1)$, otherwise
$F_x^2(x_1,y_1)=D_2(y_1)=0$, $x_1=-{B_1(y_1)/ 2B_2(y_1)}$. Thus,
for any $y$-vertex $(x_1,y_1)$, the number $y_1$ is a root of the
polynomial $D_2(y_1)$. Conversely, from the equalities \re{F_xD2}
and $y'F_y=-F_x$ one can observe that for each such root $y_1$ of $D_2$,
the point $(x_1,y_1)$ is either a $y$-vertex (i.e. an extreme point in
the $y$-direction) or a singular point. But as
we saw in subsection 
\ref{BiFol}, a singular point can occur only if a root of $D_2$ is
multiple, i.e. $\Delta=g_2^3-27 g_3^2=0$; this is an exceptional
case. Therefore, the following can be considered as proven:
\begin{pr}\label{ExtRoots}
For a generic case only the roots of $D_2(y_1)$ and nothing else can be
the $y$-coordinates of $y$-vertexes of the curve. Likewise, for the
generic case only the roots of $D_1(x_1)$ and nothing else can be the
$x$-coordinates of $x$-vertexes.\end{pr}

Now we can provide a geometrical interpretation of the cases (i), (ii),
(iii), a proof of which can be easily obtained through analyses of expressions \re{D_1} and \re{D_1D_2}.
\begin{pr}\label{StatRoots} 0). The real curve in cases (i), (ii) and (iii)
vanishes only in the case (i) under the condition $\,\tilde b<2,\,a>0$, that is equivalent to $\,b^2<(a+1)^2,\,a>0.$

1). In case (i,) of a nonvanishing real curve, the equation $D_1(x)=0$
(and the equation $D_2(y)=0$) has no real root only in the case
$\,\tilde b<2,\,a<0$ (there is no vertex here).

2). In case (i), the equation $D_1(x)=0$ (as well as the equation
$D_2(y)=0$) has four real roots under $\,\tilde b>2$ (there are four
$x$-vertices and four $y$-vertices).

3). In the case (ii), the equation $D_1(x)=0$ (as well as the equation
$D_2(y)=0$) has two real roots (there are two $x$-vertices and two
$y$-vertices).

4). In the case (iii) for $a>0\ (\hat b>2)$, the equation
$D_1(x)=0$ has four real roots and the equation $D_2(y)=0$ has no
real roots (there are four $x$-vertices and no $y$-vertices).

5). In the case (iii) for $a<0\ (\hat b<-2)$, the equation
$D_1(x)=0$ has no real root and the equation $D_2(y)=0$ has four
real roots (there are four $y$-vertices and no $x$-vertices).
\end {pr}

If we allow values $x=\infty;\ y=\infty$, i.e. consider our
real curve on a thorus $S^1\times S^1$, then the number of
$x$-vertexes and the number of $y$-vertexes will not be changed by the
M\"obius transformations \re {proj_xy}, and, thus, one obtain

\begin {pr} The linear-fractional changes  \re {proj_xy} do not lead out the curve of its class 0)-5) of the proposition \ref{StatRoots}.
\end {pr}

\subsection{John mapping of a biquadratic curve and periodicity} \label{Jbc}

First, we will show how the complex John algorithm works for the curve
$C$. Let us consider a case of some more general curve.

Let $\tilde C\subset \Bbb C^2$ be a complex curve described
parametrically via \be x(t) = \phi(t),\quad y(t) = \phi(t+\ve),\ \
t\in\Bbb C \lab{xy_phi} \ee where $\ve$ is a nonzero complex
parameter. Assume that $\phi(z)$ is an even periodic meromorphic
function, i.e. $ \phi(-z) = \phi(z), \ \phi(z+T) =\phi(z)
$ with some (complex) constant $T$. And let $C=\tilde C\cap \Bbb
R^2$ be a real curve given by means of a contraction of the
functions $x,y$ on some set $S\in\Bbb C$.


We will assume that the curve $\tilde C$ satisfies the condition:
\begin{eqnarray}\label{CDom}&  \ \mbox {\it The curve $\tilde C$
is nondegenerate and every straight line of the form $x=x_0$\qquad
}\nonumber \\  &\mbox{\it 
or $y=y_0$ intersects the curve $\tilde C$ at not more than two
points
\qquad\qquad\qquad}\end{eqnarray}

Let us perform the substitution  $t \to -t$. Now, considering the fact that $\phi(t)$ is even,
we can notice that $x(-t)=x(t)$, but, in general, $y(-t) \ne y(t)$. This
means that the point $(x(t),y(-t))$ on the curve is obtained as
the second intersection point of $C$ with the "vertical" line
passing through the initial point $(x,y)$. Thus, the transformation
$t \to -t$ is equivalent to the involution $I_1$ in the John
algorithm. Quite similarly it can be said that the transformation $t \to -t - 2\ve$
leaves the coordinate $y$ on the curve $C$ intact, whereas  $x$
is transformed to another point, which is located at the intersection of the curve $C$ with "horizontal" line. Thus, the transformation $t \to -t - 2\ve$
is equivalent to second involution $I_2$ in the John algorithm.
Obviously, for the John mapping $T$, we have the following \be T=I_2I_1
\leftrightarrow t \to t - 2 \ve, \quad  T^{-1}=I_1I_2
\leftrightarrow t \to t + 2 \ve, \lab{John_t} \ee Thus, we proved
the following
\begin{pr}\label{meromJohn}
For curve $\tilde C$, given by \re{xy_phi}, with an even periodic
meromorphic function $\phi$ satisfying the condition \re{CDom}, the
complex John mapping $T$ is equivalent to a shift of the parameter
$t$ by the step $-2 \ve$.
\end{pr}
Now the periodicity condition of the complex John mapping takes the form: \be 2n\epsilon =mT, \lab{PerCon1}\ee where $n,m$ are some integers.

In the case of the biquadratic curve, there are two periods
$2\omega_1$, $2\omega_2$, and in this case, the periodicity condition
takes the form: \be n\epsilon =m_1\omega_1+m_2\omega_2,
\lab{PerConTwo}\ee where $n,m_1,m_2$ are some integers.

Now we have only two cases which can be stated as follows: either all points of $C$ are not
periodic (for the John mapping) or every point is periodic having the
same period (which is equal to $n$, if either the condition \re{PerCon1}
or \re{PerConTwo} is satisfied).


Let us consider the last canonical form of the previous subsection -- the
Euler-Baxter curve \re{E_Bax} that is parameterized through
\re{Bax_par}. It is important that for the given $x=\sqrt{k} \,
sn(t;k)$, there be two values of $y$ corresponding to two values
of $\eta$: $y_1 = \sqrt{k} \: sn(t+\eta;k)$
 and $y_2 =
\sqrt{k}\: sn(t-\eta;k)$. These two values $y_{1,2}$ correspond to
two intersection points of the line $x=x_0$ with the curve
\re{E_Bax}.

Thus, all points $M_n$ of the complex John algorithm can be found through the following:
\be M_n =(\sqrt{k}\: sn(t+2 \eta n -3 \eta ;k),\ \sqrt{k}\: sn(t+2
\eta n -2 \eta ;k)) \lab{M_n_Bax} \ee

The periodicity condition with the period $n$ will be \be 2 \eta \,n =
4Km_1 + 2iK'm_2 \lab{per_eta} \ee where $m_1,m_2$ are some integers, and $K,K'$ are constants defined through $K=K(k)=\int\limits_0^1
{dt}/{\sqrt{(1-k^2t^2)(1-t^2)}}$, \ \ \ $K'=K(k'),\ k'^2=1-k^2$, where the notations are well-known.

For writing similar expressions for the real case, real parameterizations must be written down by means of real-valued elliptic functions.

First, let us consider  {\bf the case $a>0,\ c=+1$ of real
$x$-$y-$symmetric curve} \re{E_Bax} and use the parameterization
\re {Bax_par}. Above, in the proposition \ref{StatRoots} (case 0),
it was noted that under $a>0$, the condition $b^2>(a+1)^2$ (or $\tilde
b>2$) would be necessary and sufficient for existence of the real curve
\re{Bax_par}. Therefore, by assuming that $b^2>(a+1)^2$, it follows that $k + k^{-1}\ge
2$, and, hence, it is possible to choose $k=k_1,$ $0<k_1<1$. From the second equality \re{k_et} one obtains that $sn(\eta,k)$ is pure imaginary. This
means that $\eta = 2mK+\theta i$ with an integer $m$ and a real $\theta$
(see \cite{BatErd}). Moreover, $x$ must be real, which implies that either
$t=nK'i+\tau$ or $t=(2n+1)K+i\tau$ with an integer $n$ and real
$\tau$ (\cite{BatErd}). But $y$ must be also real, i.e.
$t+\eta=K'n_1i+\tau_1$ or $t+\eta=(2n_1+1)K+i\tau_1$ with an
integer $n_1$ and a real $\tau_1$. Now we will add expressions for $t$
and $\eta$, and show that only a variant with the second formula for $t$ and $t+\eta$: $\tau_1=\theta+\tau$ is possible. Thus,
we assume that the parameter $t$ isn't real: $t=\pm
K+i\tau$; however, the values $x,y$ in the paramentization
\re{Bax_par} are real and the parameter $\tau$ is also real, here the
sign before $K$ determines a branch of the curve, which has a shape of a two bounded oval in this case.
The relations \re{k_et} have a
denumerable set of solutions $\eta$, and we choose $\theta=$ Im $\eta$
as a minimal positive number of all Im $\eta$.
Now the periodicity condition \re{per_eta} can be written in the
form 
($m=m_2$): \be \frac{\theta}{K'}=\frac{m}{n}
\in\Q\,\lab{per_theta} \ee and from relations \re{k_et}, one
obtains that the number $\theta$ is a minimal positive solution of
the equation $\mbox{sc}(\theta,k')=1/\sqrt{ak}$ \quad (here and
below the following adoptions, as usual, will be made $sc=sn/cn,\ ns=1/sn$, and so on).

Hence, we proved the foolowing

\begin{pr}\label{thetaJohn1} The condition
\re{per_theta} is a criterion of periodicity of the John mapping
$T$ in the case $a>0,\ c=+1$ of real $x$-$y-$symmetric curve
\re{E_Bax}. In this case, all points of our curve are periodic with the
same minimal period $n$ in a nonreducible fraction $\frac{m}{n}$.
The number $m$ is a number of full turns of the curve that the
mapping $T^n$ performs.
\end{pr}

Let us also {\bf another cases of real $x$-$y-$symmetric curves}
\re{E_Bax}. Let us take an advantage of a list of cases given in the work
\cite{IatrouRob} (note that the latter work contains errors for the
above case $c=1,\,a>0$, so here we used formulae from \cite{Bax} instead). Now let us write down the final parameterization formulas, in which the modulus $k$ is given by $\hat b$ (remind that $k'\,^2+k^2=1$) and the shift $\eta $ is given by $a$ or $b$ (as they are related). Here $\hat b=(b^2-a^2-c)/a$.
\begin{equation} \begin{aligned}&\text{Case } a>0,\,c=-1,\, :\quad  x=\sqrt{k/k'} \, cn(t;k),
\quad y=\sqrt{k/k'} \: cn(t\pm\eta;k),\\
&\text{where }\ k/k'-k'/k=\hat b,\ \ a=ds^2(\eta,k)/kk',\
b=-cs(\eta,k)\,ns(\eta,k)/kk'.\quad\end{aligned} \lab{par13}
\end{equation}

\begin{equation} \begin{aligned}\text{Case  } a<0,\ \ &c=-1,\,:\quad  x=\sqrt{k/k'} \, nc(t;k),
\quad y=\sqrt{k/k'} \: nc(t\pm\eta;k),\\
\text{where }\ k'/k-&k/k'=\hat b,\quad a=-ds^2(\eta,k)/kk',\
b=cs(\eta,k)\,ns(\eta,k)/kk'.\quad\end{aligned} \lab{par58}
\end{equation}

\begin{equation} \begin{aligned}\text{Case } a<0,\,c= 1,\,\hat b<-2&:\   x=\sqrt{1/k'} \, cs(t;k),
\ y=\sqrt{1/k'} \: cs(t\pm\eta;k),\\
\text{where }\  -1/k'-k'=\hat b,\ \ &a=-cs^2(\eta,k)/k',\
b=ds(\eta,k)\,ns(\eta,k)/k'.\quad\end{aligned} \lab{par4}
\end{equation}

\begin{equation} \begin{aligned}\text{Case } a<0,\,c=1,\,\hat b>2: \hspace{8cm}
\\
\text{unbounded part } x=\sqrt{1/k} \ ns(t;k),
\ \ y=\sqrt{1/k} \ ns(t\pm\eta;k),
\hspace{1.1cm}\\
\text{bounded part } x=\sqrt{k} \ sn(t;k),
\quad y=\sqrt{k} \: sn(t\pm\eta;k),\hspace{2.2cm}
\end{aligned}
\lab{par7}
\end{equation}         
\quad where 
$1/k+k=\hat b,\quad a=-ns^2(\eta,k)/k,\
b=cs(\eta,k)\,ds(\eta,k)/k.$ 

\vskip 5pt The separation present in the last case is related to the fact that the curve has one bounded oval branch, and likewise it contains some other branches, which are unbounded. The last four cases were also given in the work \cite{IatrouRob}.

The remaining cases were considered in the work \cite{IatrouArx}:
The case  $a<0,\,c=1,\,|\hat b|<2$
of real symmetric curves \re{E_Bax} 
can be transformed to the case \re{par4} via the following substitution
$$(x,y)= \left({{1-\bar x}\over{1+\bar x}},{{1-\bar y}\over{1+\bar
y}}\right). \text {\ Then,\ } \bar a= {{1+a-b}\over {1+a+b}},\
\bar b={{2(1-a)}\over{1+a+b}},\ \bar{\hat
b}=2+{{16a}\over{b^2-(a+1)^2}}.$$

\begin{equation} \begin{aligned}\text{Case } a<0,\,c= 1,\,|\hat b|<2&:\   \bar x=\sqrt{k'} \, sc(t;k),
\ \bar y=\sqrt{1/k'} \: sc(t\pm\eta;k),\\
\text{where }\  -1/k'-k'=\bar{\hat b},\ \ &\bar
a=-cs^2(\eta,k)/k',\ \bar
b=ds(\eta,k)\,ns(\eta,k)/k'.\quad\end{aligned} \lab{par asym}
\end{equation}


\vskip 5pt Let us consider cases of {\bf real $x$-$y-$asymmetric curves \re{2_2F}, the
case (iii)} (see the statement \ref{StatReduce}):
\begin{equation} \begin{aligned}\text{Case } a>0,\,c= 1,\,(\hat b>2)&:\   x=\sqrt{k'} \, nd(t;k),
\ y=\pm\sqrt{1/k'} \ dn( t\pm\eta;k),\\
\text{where }\  1/k+k=\hat b,\ \ &a=k'nd^2(\eta,k),\ b=k^2
sd(\eta,k)\,cd(\eta,k).\quad\end{aligned} \lab{par a1}
\end{equation}
\begin{equation} \begin{aligned}\text{Case } a<0,\,c= 1,\,(\hat b<-2)&:\   x=\sqrt{k'} \, sc(t;k),
\ y=\sqrt{1/k'} \: cs(\pm t+\eta;k),\\
\text{where }\  -1/k-k=\hat b,\ \ &a=-k'\,nd^2(\eta,k),\ b=k^2
sd(\eta,k)\,cd(\eta,k).\quad\end{aligned} \lab{par a2}
\end{equation}

All of these cases, i.e. \re{par13}-\re{par a2}, allow calculations, similar to those described above, and the calculations finally led to the formula  \re{per_theta}. Such calculations were completed by our collaborators
M.V. Belogljadov and A.A. Telitsyna (\cite{Bel },\cite{Telitsina}), and their work produced the following results:

a common answer for cases 
\re{par13}, \re{par4}, \re{par7}, \re{par asym}, \re{par a2}


\begin{equation}
{\mbox{Re }\eta\over{2K}}={m\over n}\in \Bbb Q; \lab{res1}%
\end{equation}

an answer for the case \re{par58} in the form
\begin{equation}
{\mbox{Re }\eta\over{K}}={m\over n}\in \Bbb Q; \lab{res2}%
\end{equation}

and, at last, for the case \re{par a1} an answer was given in the form:
\begin{equation}
{\mbox{Im }\eta\over{2K'}}-1/2={m\over n}\in \Bbb Q. \lab{res3}%
\end{equation}

Here $k$ and $\eta$ are calculated from relations
\re{par13}-\re{par a2}, depending on a case to be involved
Thus, the following proposition can be made:

\begin{pr}\label{etaJohn}
The conditions 
\re{res1}-\re{res3} are criteria of periodicity of the John
mapping $T$ in the corresponding real cases \re{par13}-\re{par a2}.
The number $n$ is a period of dynamical system and $m$ is the
number of full turns of mapping $T^n$.
\end{pr}

Note that correspondence to the cases of proposition \ref{StatRoots} can be easily observed, and it must be noted that such presence of several cases corresponding to a case of proposition \ref{StatRoots} is related to technical reasons (namely, to the choice of the infinite point on the projective line).

\section{The Poncelet problem}
\setcounter{equation}{0}

In this section we demonstrate a nice correspondence between the
mapping by F.John and the famous Poncelet problem (the Poncelet
porism) for two conics. We will start by recalling the Poncelet porism
in a well-known form.

\subsection{The Poncelet porism in the form of two circles
}\label{HisRemPP}

To begin with, let a circle $A$ lie inside another circle $B$.
From any point on $B$, let us draw a tangent to $A$ and extend it to $B$.
From the intersection point, let us draw another tangent, etc. For $n$
tangents, the result is called an $n$-sided Poncelet transverse.
This Poncelet transverse can be closed for one point of origin,
i.e. there exists one circuminscribed $n$-gon (which is inscribed
into the outer circle and, at the same time, circumscribed around
the inner circle). We may begin with a polygon, which will be
understood as a set of straight lines bridging, in a sequential
manner, a given cyclic sequence of points (i.e. vertices) on the
plane. If there exist two circles,
inscribed and circumscribed for this polygon, then this polygon is
called a bicentric polygon. Note that sides of the polygon are
allowed to intersect, and the intersection point will not
necessarily be a vertex. Furthermore, the inscribed circle does not
necessarily contact a segment located between vertices, and the contact
point can lie on an extension of the side and therefore the circles
may intersect. Bicentric polygons serve as popular objects of
investigations in geometry. This is the most well-known form
of the Poncelet porism.

If we denote a radius of the inscribed circle as $r$, a radius of
the circumscribed circle as $R$ and a distance between the circumcenter
and incenter for a bicentric polygon as $d$, then these three numbers
will not be arbitrary and along with $n$, they will have to satisfy
certain relations. So, for the case of a triangle, such relation is
sometimes called as the Euler triangle formula $R^2-2Rr-d^2=0$,
which was well-known even in Babylon for some particular cases.
One of popular notations for such relations (which is necessary and
sufficient for existence of a bicentric polygon) can be given in terms
of additional quantities
$$a=\frac{1}{R+d}\ ,\quad b=\frac{1}{R-d}\ ,\quad c=\frac{1}{r}\ .$$
So, for a triangle above, the Euler formula has the form: \
$a+b=c$, for a bicentric quadrilateral, radii and a distance are
related to each other by the equation\quad $a^2+b^2=c^2$.
The relationship for a bicentric pentagon is \quad $4(a^3+b^3+c^3)=(a+b+c)^3$.
For a general case, one can introduce the following numbers
$$\lambda=1+\frac {2c^2(a^2-b^2)}{a^2(b^2-c^2)}\
,\quad \omega=\cosh^{-1}\lambda\ ,\quad k^2=1-e^{-2\omega},\quad
K=K(k)\ (\mbox{see s.\ } \ref{Jbc})$$ and then the relationship
can be written through elliptic functions in the form
\be\mbox{sc}\left(\frac{K}{n},k\right) =\frac{c\sqrt
{b^2-a^2}+b\sqrt{c^2-a^2}}{a(b+c)}\lab{bicentr} \ee (Richelot
(1830) was the first to give a criterion, which was difficult
and imperfect, then another criterions were published by different
authors and, at last, the written criterion by S.M. Kerawala was
published in his work \cite{Kerawala} in 1947 year).

%
%
%

\subsection{Setting of the Poncelet problem}

Let us recall, for simplicity, the Poncelet problem \cite{Berger} for the case of two ellipses, as it was introduced by Jean-Victor Poncelet himself.
We will take two arbitrary 
ellipses $A$ and $B$, with $A$ located inside $B$. Let us have an arbitrary
point $Q_1$ on the ellipse $A$ and draw
a tangent straight line to $A$ at the point $Q_1$. This tangent
crosses the ellipse $B$ at two points $P_1$ and $P_2$, , and, besides, we will assume that it crosses $P_1$
before $P_2$ with respect to the standard orientation. Then let us take
the point $P_2$ on $B$ and draw the second tangent to the ellipse
$A$. We denote as $Q_2$ the point on $A$, in which the tangent
contacts $A$. This tangent crosses the ellipse $B$ in two
points $P_2$ and $P_3$. We will take the point $P_3$ and repeat the
procedure. Then we obtain a mapping $U_B:B\to B$, which acts in compliance with the following rule $U_B:P_k\to P_{k+1}$, which will be called below as the Poncelet mapping. Moreover, we obtain the
mapping $U_A:A\to A$ that acts by the rule $U_A:Q_k\to Q_{k+1}.$ More
precisely, due to the fact that the definition of $U_B$ uses a point $Q\in
A$, we can introduce two mappings $I_A,I_B: 
\tilde C\to \tilde C$,\ $\ \tilde C:=\{(Q,P)\in A\times B\ |\ P
\mbox{ lies on a tangent line to } A \mbox{ at } Q \}$, working by
the rules $I_A:(Q_1,P_2)\to (Q_2,P_2)$, $I_B:(Q_1,P_1)\to
(Q_1,P_2)$. The mappings $I_B,I_A$ generate a composition
$U=I_B\circ I_A$, which is similar to the John mapping $T$. We also
obtain two sets of points $P_n$ and $Q_n$ on the ellipses $B$
and $A$, respectively.

The mapping $U_B$ has an inverse and generates a discrete
dynamical system or, in other words, an action of the group $\Bbb
Z$ on $B$, as it was above for the John mapping. An orbit of this action
is a set of points $P_k,\ k=...,-1,0,1,2,...$ and
$P_k=U_B^{k-1}P_1$. Now note that in a general case of a
disposition, the ellipses can be intersecting, and we must start
from a point $Q_1$ on the ellipse $A$, which is located inside $B$ and
we can determine the mappings $U_B$ and $U_A$ in the same way. In this case we can encounter a situation, in which the tangent straight line intersects the ellipse $B$ only at one point, and in that case, this point will be considered as a double point. The ellipses can be tangent in one or two points, or they can
be tangent and intersecting at the same time, or they can be
nonintersecting and lie one outside other, or finally, they can be
arbitrary irreducible conics.

The approach is the same in all these cases. Note that one can consider a case, in which the conic $B$ is reducible, for example, in the situation of two non-tangential different straight lines. Also note that every projective transformation of the plane transforms a Poncelet mapping of conics into the same mapping of their images; therefore, we can restrict ourselves to a case, in which one of the conics is a unit circle.

The first interesting problem to solve here would be to find a way of describing these crossing points explicitly. The problem was solved by Jacobi and Chasles, who showed that the sequences $P_n$ and $Q_n$ can be parameterized by means of elliptic functions. The second problem, the so-called Poncelet porism or the big Poncelet theorem (see \cite{Berger}), is related to proving that if a particular trajectory of action on conics is a closed trajectory (i.e., if $P_N = P_0$ for some $N > 2$), then this property does not depend on a choice of an initial point $Q_1$ on the conic $A$ (or on a choice of $P_0$). A modern treatment of this problem from an algebro-geometric point of view can be found, for example, in \cite{GH}. We follow a different approach. If we introduce standard rational parameterizations of our conics (see below \re{par_D}, \re{par_C}), then the parameters $x$ and $y$ of the points $Q_1, P_1$ can be proven to be related via a polynomial equation $F(x, y) = 0$. In a generic situation, $F(x, y)$ is a polynomial of the exact order two with respect to each of the variables. Indeed, for a nondegenerate situation, a tangent to the conic $A$ at the point $x$ must intersect $B$ in two distinct points, and, moreover, from the point $y$ on $B$, there exist two distinct tangents to $A$. Conversely, for any polynomial $F(x,y)$ of the order two on $x$ and $y$, it is possible to find two conics $A$ and $B$ parameterized through \re{par_D} and \re{par_C}. These arguments are sufficient for building trajectories of the John dynamical system. Note that similar considerations were exploited in \cite{FrRag} for showing that the tangent and the intersection points belong to an elliptic curve. The authors of \cite{FrRag} also introduced two involutions $I_1$ and $I_2$, which are equivalent to our involutions in the John $T$-algorithm for the Euler curve. In our approach, we derive the curve $F(x,y) = 0 $explicitly and then we conduct its investigation. One can note connections to Gelfand's question \cite{King} and the elastic billiard  \cite{Ves},
\cite {Laz}.

\subsection{Transition to the John mapping on a biquadratic curve}\label{PonJohn}
Let us find an explicit expression for the Poncelet
mapping $U_B$. We will introduce the standard rational parameterization
of an arbitrary conic \cite{Berger} that can be found even in the case
when the conic is reducible. Assume that the conic $A$ is described by the
coordinates $\xi_0, \xi_1,\xi_2$ of a two-dimensional projective
plane, expressed via
\be\sum
_{i,j=0}^{2} \widetilde A^{ij}\xi_i\xi_j=0 \lab{A^}
\ee and, as is well-known, the conic is irreducible, if the
matrix $\widetilde A$ is nondegenerate.

Note that we will use lower indices for vectors (contravariant tensors), and upper indices for low used covectors (covariant tensors), since indices are used often for power exponents.

Corresponding affine coordinates will be denoted as $\xi,\eta$:
$[\xi_0:\xi_1:\xi_2]=[1:\xi:\eta]$, and it is implied that $\xi_0\ne 0.$
Then it is possible to find polynomials $E_0(x),E_1(x),E_2(x)$
with $\deg(E_i(x)) \le 2$ such that \be \xi_i=E_i(x),\quad\mbox{\
or\ }\quad\xi=\frac{E_1(x)}{E_0(x)},
\quad\eta=\frac{E_2(x)}{E_0(x)}. \lab{par_D} \ee Quite
similarly, the conic $B$ can be parameterized as \be
\xi_i=G_i(y),\quad\mbox{\ or\ }\quad
\xi=\frac{G_1(y)}{G_0(y)},\quad \eta=\frac{G_2(y)}{G_0(y)},
\lab{par_C} \ee where $G_i(y)$ are some other polynomials of
degrees, not exceeding two. Thus the value of parameter  $x$
completely characterizes a point on the conic $A$, and the value
of $y$ completely characterizes a point on the conic $B$.

Now let us consider a more general case.

\begin{lem}Let the curves $A$ and $B$ be given parametrically by \re{par_D},
\re {par_C} with some smooth functions $E_i,G_i$ and
$E_0\not\equiv 0, G_0\not\equiv 0$ on each interval. The point
$P=G(y)$ lies on a tangent line to $A$ at $Q=E(x)$ iff
$F(x,y)=0$ where \be F(x,y)=\left|\begin{array}{ccc}
 E_0(x) & E_1(x) & E_2(x)\\
E'_0(x) & E'_1 (x) & E'_2(x) \\ G_0(y) & G_1(y) & G_2(y)
\end{array} \right |\ .\lab{Fdet}\ee
\end {lem}

\begin{proof}
 The affine tangent line
$L$ to the curve $A$ at point $Q$ with parameter $x$ has the
direction vector $\tau=(\frac{d \xi}{d x},
\frac{d \eta}{d x})
$. Therefore, the affine point $P$ satisfies the equality
$$\overrightarrow {OP}(y)- \overrightarrow {OQ}(x)=k\tau$$ with some
$k\in\Bbb R$ and an origin $O$. The last equality is equivalent to
the complanarity condition of vectors $\vec E=(E_0(x),E_1(x),E_2(x))$, $\vec
E'=(E'_0(x),E'_1(x),E'_2(x))$ and $\vec G=(G_0(y),$
$G_1(y),G_2(y))$ in the bundle space $\Bbb R^3\setminus\{0\}$ of
the projective fiber bundle $\Bbb R^3\setminus\{0\}\to\Bbb
{RP}^2$, where the prime denotes derivative with respect to
$x$. Indeed, collinearity of the vectors $\overrightarrow
{OP}(y)- \overrightarrow {OQ}(x)$ and $\tau$ means that
$$
0 = \left |
\begin{array}{cc} \frac{E_1(x)}{E_0(x)}-\frac{G_1(y)}{G_0(y)}
& \frac{E_2(x)}{E_0(x)}-
\frac{G_2(y)}{G_0(y)}  \\
\left(\frac{E_1}{E_0}\right)^\prime (x) &
\left(\frac{E_2}{E_0}\right)^\prime (x)
\end{array} \right |= \left|
\begin{array}{ccc} 1 &
\frac{E_1}{E_0} & \frac{E_2}{E_0} \\
0 & \frac{E_1}{E_0}-\frac{G_1}{G_0} & \frac{E_2}{E_0}-
\frac{G_2}{G_0}
\\ 0& \frac{E'_1}{E_0}-\frac{E_1E'_0}{E_0^2} &
\frac{E'_2}{E_0}-\frac{E_2E'_0}{E_0^2}  \end{array} \right |=
$$
$$
=-\left|
\begin{array}{ccc}
1 &
{E_1}/{E_0} & 
{E_2}/{E_0} \\
1 &
{G_1}/{G_0} &  
{G_2}/{G_0}
\\ {E'_0}/{E_0} &
{E'_1}/{E_0} & 
{E'_2}/{E_0}
 \end{array} \right |= \left|
\begin{array}{ccc}
E_0(x) & E_1(x) & E_2(x)\\
E'_0(x) & E'_1 (x) & E'_2(x) \\ G_0(y) & G_1(y) & G_2(y)
\end{array} \right |\frac{1}{G_0(y)E_0(x)^2}\ .
$$
\end{proof}

Let us return to our case of conics. One can see $F(x,y)$ has the form \be
F(x,y)= M_0(x) G_0(y) + M_1(x) G_1(y) + M_2(x) G_2(y)
\lab{Z_Ponc2} \ee with polynomials $M_i(y)$ defined as \be M_i(x)
= \epsilon_{ikl} (E_k'(x)E_l(x)- E_k(x)E_l'(x)), \quad
i,k,l=0,1,2, \lab{M} \ee where $\epsilon_{ikl}$ is a completely
antisymmetric tensor. One can easily check that $\deg(M_i(x)) $
$\le 2$, and therefore the curve $F(x,y)=0$ is a biquadratic curve of
the form \re{2_2F}.

Note that the equality \re{Z_Ponc2} can be written as
$F(x,y)=(\,\vec M(x),\vec G(y)\,)$ with the scalar product
$(\cdot,\cdot)$ or $F=\langle\, 
\vec M(x), \vec G(y)\,\rangle$ with the pairing $\langle
\cdot,\cdot\rangle$ (in this case $\vec G\in \Bbb R^3, \vec M\in
(\Bbb R^3)^*$\ ) and
$F=(\vec E(x),\vec E'(x),\vec G(y))$ with the mixed 
product in $\Bbb R^3$. Further, we introduce vectors $\vec x=colon
(1,x,x^2), \vec y=colon (1,y,y^2)$ and matrixes $ E,G$ by the
rules $E_i(x)=\sum_{j=0}^2 E_{ij}x^j=( E\vec x)_i$, $\vec G
= G\vec y$, and then the decomposition \re{Z_Ponc2} can be written in
the form $F(x,y)=( M\vec x, G\vec y)$, where $\vec M(x)=M\vec x$.
This implies that $F=(\vec x, M^* G\vec y)=( G^* M\vec x,\vec y)$.
Comparing with the decomposition in the forms \re{A_form} and
\re{B_form}, we obtain \be A= G^* M,\ B= M^*G=A^*, \lab{ABGM}\ee
where the matrix $A$ is obtained in the same way as the matrix $E$
above; moreover $A=(a_{ik})$ with the matrix from \re{2_2F}.

In a case of irreducible conics the matrix $A$ is nondegenerate.
Indeed, if the matrix $A$ is degenerate then either the matrix $M$
or the matrix $G$ is degenerate by virtue of \re{ABGM}. The
degeneracy of $G$ means a linear dependence of polynomials
$G_i(y)$, that is, the conic $B$ will be a straight line in this
case. The degeneracy of $M$ means that there exist constants
$c_0,c_1,c_2$ such that
$$ \left|\begin{array}{ccc}
 E_0(x) & E_1(x) & E_2(x)\\
E'_0(x) & E'_1 (x) & E'_2(x) \\ c_0 & c_1 & c_2
\end{array} \right |\equiv 0,\ $$ i.e. a linear
dependence of polynomials $E_i(x)$ and the conic $A$ will be a
straight line.

Let us return to the Poncelet construction. We have obtained the
parameters $x_1$ of the point $Q_1$ and $y_1$ of the point $P_1$,
satisfying the equation $F(x_1,y_1)=0,$ with $F$ from \re{Fdet}. Note
now that instead the point $P_1$, we could write the point $P_2$
with the parameter $y_2$ in the Poncelet construction and have the
same equation $F(x_1,y_2)=0.$ We obtain the first result: for any
point given by a parameter value $x_1$ on the conic $A$, the
points with parameters $y_1$ and $y_2$ of the intersection points
of the tangent line $L_1$ at $x_1$ with $B$ are determined as two
roots of the quadratic equation: $ F(x_1,y)=0.$
Thus, by identifying a point with its parameter value, we can say that the
Poncelet mapping $U_B$ maps the point $y_1$ into the point $y_2$
and, hence, it (more precisely, the mapping $I_B$) coincides with
the John mapping $I_1$
from section 2. 
Similarly, the mapping $U_A$ maps the point $x_1$ into the
point $x_2$ and, hence, it (the mapping $I_A$) coincides with the
John mapping $I_2$, and the mapping $U$ coincides with the John
mapping $T$ on the curve $C$ \re{2_2F}.


We have proven the 
following


\begin{pr}\label{JohnPonc1} Each pair of distinct irreducible real conics $A$,
\re{par_D} and $B$, \re{par_C} generates a biquadratic real
polynomial $F(x,y)$ of the form \re{2_2F} with nondegenerate
matrix $A$ by means of \re{Fdet} such that the point $y\in B$ lies
on the tangent at the point $x\in A$ if and only if $F(x,y)=0.$
The Poncelet mapping $U_B$ gives us (and can be obtained from) the
mapping $U:\tilde C\to \tilde C$ that coincides with the John
mapping $T$ on the curve $C$ \re{2_2F}.
\end{pr}

Let us prove the following converse

\begin{pr}\label{JohnPonc2} For any biquadratic real polynomial $F(x,y)$ of the form \re{2_2F} with the nondegenerate matrix $A$ and for its every
decomposition \re{Z_Ponc2} with given polynomials $G_i(y),M_i(x)$
of the second order, there exists a unique projective set of
polynomials $E_i(x)$ of the second order such that relations
\re{M} (and \re{Fdet}) hold, and, hence, we can relate with any
such curve $F(x,y)=0$ and its decomposition \re{Z_Ponc2} a pair of
conics $A$ and $B$ parameterized as in
\re{par_D} and \re{par_C}. 
\end{pr}

{\it Proof}. Let us have a decomposition \re{Z_Ponc2} of a given
biquadratic polynomial $F$. Then we have a parameterization of the
conic $B$ by means of projective coordinate $G_0,G_1,G_2$ and our first aim
is to find polynomials $E_0,E_1,E_2$ such that the
equalities \re{M} will hold. In other words, we must find a polynomial
solution of the following system of ordinary differential equations \be
E\times E'=M\lab{rr}\ee with a known vector $M=(M_0(x), M_1(x),
M_2(x))$, an unknown vector $E=(E_0(x),$ $ E_1(x), E_2(x))$ and
the vektor product $\times$ in $\Bbb R^3.$

One can think of this equation as of a problem to find a curve (more precisely, a tangent vector of a curve), if its binormal is known.

Now we will need
the following

\begin{lem} Consider the differential equation in $\Bbb R^3$
\be r\times r'=b\lab{rrl}\ee with a known smooth vector function
$b=(b_1(t), b_2(t), b_3(t))$ and unknown vector function
$r=(r_1(t),$ $ r_2(t), r_3(t))$ depending on a real parameter $t$.
If $(b,b',b'')<0$, then the equation \re{rrl} does not have any smooth
solution. If $(b,b',b'')>0$, then a solution of the equation
\re{rrl} exists and it can be expression only as $r=\pm(b,b',b'')^{-1/2}\,b\times b'$. 
If $(b,b',b'')= 0$ on an interval, then only four cases are
possible with some parameter change $s=s(t)$: 1) $r\equiv
r_0=const$, 2) $r= r_0+r_1s,$ $\ r_0,r_1$ are constant vectors. 3)
$r=r_0 e^s+r_1e^{-s}$ and 4) $r=r_0\cos s+r_1\sin s$ where $r_0, r_1$ are the same as before.
\end{lem}

\begin{proof} First, note that the Jacobi determinant of the system \re{rrl}
of ordinary differential equations with respect to $r'$
identically equals zero, so a standard theory of
differential equations systems does not work for this system.

I). Let 
$b\times b'\not \equiv 0$ on an interval. We observe that the
solution must be only of the form $r=v(t)\,b\times b'$ because
from \re{rrl}, one can easily obtain the following: $b\cdot r=0,\ b\cdot r'=0$, so
$b'\cdot r=0.$

The scalar $v$ is still unknown.

After substitution of such
$r$ in \re{rrl}, we obtain $v^{2}(b,b',b'')b=b$ so that there
exists no solution of \re{rrl}, when $(b,b',b'')< 0$. The equality
$(b,b',b'')= 0$ on an interval gives $b=0$, which is a contradiction.

II). Now let $b\times b'\equiv 0$ on an interval. If $b\not \equiv
0$ and $b'\not \equiv 0$, then $b(t)=\mu(t)b_0$ with $\mu$ being a scalar
function, and $b_0$ being a constant vector, which means that the curve
$r$ is a plane curve. Note that we simultaneously examine the case
$b\not \equiv 0$ and $b'\equiv 0\ (\mu\equiv 1)$. Let us substitute this
$b$ into \re{rrl} and choose another parameter $\tilde t$ such that
$r\times r'=b_0$. It follows from this, that $r\times r''=0$, so
$r''=\nu(\tilde t)r$ with $\nu$ being a scalar function. After a reparameterization $\tilde t\to s$ we obtain $r''=\pm r$, i.e.
$r=r_0e^s+r_1e^{-s}$ or $r=r_0\cos s+r_1\sin s$ with some constant
vectors $r_0, r_1$ of a plane, orthogonal to $b_0$. For
$r_0\ne r_1$, these solutions will describe cases of a hyperbolic rotation and an elliptic rotation of the plane, pespectively; besides, also there exist cases of a movement along a straight line (cases $b=0$
and $r_0=r_1\not=0$) and of a stationary state.
\end{proof}

{\it Continuation of the proposition proof.} We will apply the
lemma to the case of polynomial vectors of the second order $b=M$
and $r=E$. It is easy to see that for such $M$, the equality
$M\times M'\equiv 0$ implies only trivial cases 1) or 2) of the
lemma, so that $M\times M'\not\equiv 0$ and $v
^2=(M,M',M'')$. 
The latter mixed product does not depend on $x$.
Indeed,
$$
{d\over {dx}}(M,M',M'')=(M',M',M'')+(M,M'',M'')+(M,M',M''')=0
$$
as the vector $M$ consists of quadratic polynomials. Note that we
can change a sign of the mixed product $(M,M',M'')$ by changing of
the sign of $M$
. Thus, for given $M_i$ of the decomposition \re{Z_Ponc2} the
scalar $v$ exists there, and components of the derived vector
$E=vM\times M'$ are polynomials of the second order (and therefore
the mixed product $(E,E',E'')$ does not depend on $x$). $\Box$

\subsection{Projective invariance of the biquadratic curve }
One can consider the vector $E\times E'$ as a covector $E^*$ of a
dual space $(\Bbb R^3)^*$ with
components $E^j(x)=\sum
_{i=0}^{2} \widetilde A^{ij}E_i(x)$, where the matrix $\widetilde
A$ is from \re{A^}, or more precisely, $E^*$ is proportional to
$E\times E'$. Indeed, by definition, we, at once, obtain $\langle
E^*(x),E(x)\rangle\equiv 0$ and $\langle E^*,E'\rangle\equiv 0$.
The covector field $E^*(x)$ describes a field of tangent lines to $A$, i.e., a conic in the dual projective space (see, e.g., \cite{GH1}):
\be\sum
_{i,j=0}^{2} \widetilde A_{ij}E^iE^j=0 \lab{A_}\ee with an inverse
matrix $(\widetilde A_{ij})=(\widetilde A^{ij})^{-1}$. This viewpoint
gives the following representation \be F(x,y)=\langle
E^*(x),G(y)\rangle\ee and helps to understand the reason why the equality
$E\times E'=M$ with the quadratic polynomial vector $E$ implies
$M\times M'=E$ in the projective sense of $E^{**}=E$, but it
does not provide help for the case, in which the curve $A$ has an order greater than two. 

\begin{pr} For every pair of different irreducible real conics $A$,
\re{par_D} and $B$, \re{par_C}, 
and for any projective transformation $L$ of the projective plane, the
images $LA$ and $LB$ of the conics $A$ and $B$ generate the same
equation $F(x,y)=0$ and the curve \re{2_2F}. Two different pairs
of conics with the same 
curve are connected by means of a projective transformation $A\to
P_LA,\ B\to P_LB$ 
and give different decompositions \re{Z_Ponc2}.
\end{pr}

\begin{proof} In order to derive the first result, let us recall that an
arbitrary nondegenerate projective transformation $P_L$ of the
affine plane can be obtained from a nondegenerate linear
transformation $L$ of the bundle space $\Bbb R^3\setminus\{0\}$ of
the projective fiber bundle. Such linear transformation $L$ will give
a number factor  $det\,L$ for $F$ in the formula \re{Fdet}. Thus,
an arbitrary nondegenerate projective transformation of the affine
plane and the conics $A$ and $B$ will not change the equation \re{2_2F}.

Then, the given decomposition \re{Z_Ponc2} can be written in the
forms \re{A_form} and \re{B_form} and it will generate the nondegenerate
matrixes $M,G$, and then, $A,B$ my means of \re{ABGM}. Some other
decomposition 
would give matrixes $M_1,G_1$ and the same $F$ and $A$: $A= G^* M=A_1=
G_1^* M_1,$ so that the linear transformations
$L=G_1G^{-1}=(M_1^{-1})^*M^*$ and $M_1M^{-1}=(G_1^{-1})^*G^*$ would
translate the conics $A,B$ into conics $A_1=LA,B_1=LB$,
respectively, as they are projective transformations. It can be noted from \re{A^}
that under the transformation $L$ the matrix $\widetilde A$ will convert
into $L^{-1*}\widetilde AL^{-1}$, and the matrix $M$ will convert into
$L^{-1*}M$.
\end{proof}






Thus, the biquadratic curve $F(x,y)=0$ depends only on a
projective class of a pair of conics. But the conics can be present in one
of the following generic dispositions (here the order is selected to match the order used in proposition \ref{StatRoots}):

0) the conic $B$ lies inside of $A$, strictly or not, i.e. no
tangent line to $A$ cuts $B$.

1) the conic $A$ lies strictly (i.e. without a contact) inside of
$B$, so that every tangent line to $A$ cuts $B$ in two different
points.

2) the conic $A$ cuts the conic $B$ in four points and there exists
a straight line that has no common point with $A$ and $B$.

3) the conic $A$ intersects the conic $B$ at only two points
without contacts.

4) the conic $A$ lies strictly outside of $B$ and the conic $B$
lies outside of $A$.

5) the conic $A$ cuts the conic $B$ in four points and every
straight line has a common point with $A$ or $B$ (there is no
common tangent line, a hyperbola and an ellipse).








%
%
%
%
%
%
%
%
%

\vskip 5pt 
Now we would like to clarify a question how to find from $F$, what case
occurs. In order to clarify this, let us make the following
observations. If the conic $A$ intersects $B$ in a point $P$, then
from the point $P\in B$ there may be only a tangent line to $A$, so
that in the Poncelet construction there is only a point $(x,y)$
with the parameter $y$ corresponding to $P\in B$, that is, the point
$(x,y)$ is a $y$-vertex of $C$, i.e. an extreme point along a
direction of the real axis $y$. Remember that by the proposition
\ref{ExtRoots}
 {\it for generic cases
1)-5) the roots of $D_2(y_1)$ and nothing else can be the $y$-coordinates of $y$-vertices. Every such vertex corresponds to a point of
intersection $A\cap B.$} 

In a similar way, a common tangent line to conics $A$ and $B$ in the
plane $(\xi,\eta)$ gives a $x$-vertex $(x_2,y_2)$ in the plane
$(x,y)$ and we have the following: if $A_2(x_2)=0$ then $y_2=-A_0(x_2)/A_1(x_2)$
otherwise $D_1(x_2)=A_1^2(x_2) - 4 A_0(x_2) A_2(x_2)=0$,
$y_2=-{A_1(x_2)/ 2A_2(x_2)}$. And {\it for generic cases 1)-5),
real roots of $D_1(x_2)$ and nothing else can be $x$-coordinates of
$x$-vertices. Each such vertex corresponds to a common tangent
line.}

Let us recall the statement \ref{invar}, saying that the invariants
$g_2$ and $g_3$ of polynomials $D_1$ and $D_2$ are identical. Note
that the case $\Delta<0$ (i.e. $k^2<0$) corresponds only to the
above case 5) because the polynomial $D_1$ has two real roots and two
complex roots if and only if $\Delta<0$, or equivalently, if and only if
the polynomial $D_1$ is the same. We obtain the following

\begin{pr} In the case of disposition 0), each of the equations
$D_1(x_1)=0$ and $ D_2(y_1)=0$ has no real solution, and the curve
is not real. In the case 1) each of the equations $D_1(x_1)=0,
D_2(y_1)=0$ has no real solution, the curve $C$ is real and it contains
no vertex. In the case 2) each of the equations $D_1(x_1)=0$
and $D_2(y_1)=0$ has four real solutions, the curve $C$ has four
$x$-vertices and four $y$-vertices. In the case 3) each of the
equations $D_1(x_1)=0, D_2(y_1)=0$ has two real solutions, the
curve $C$ has two $x$-vertices and two $y$-vertices. In the case
4) we have four common points, so that the equation $D_1(x_1)=0$
has four real solutions, the equation $ D_2(y_1)=0$ has no real
solution, the curve $C$ has no $y$-vertex and four $x$-vertices.
In the case 5) we have four common tangent lines, so that the
equation $D_1(x_1)=0$ has no real solution, the equation $
D_2(y_1)=0$ has four real solutions, the curve $C$ has no
$x$-vertex and four $y$-vertices.

\end{pr}

Any point of contact gives a point $(x_1,y_1)$, which acts as both
$x$-vertex and $y$-vertex, i.e. $D_1(x_1)=D_2(y_1)=0$, the
conditions \re{stab} of singular point are satisfied, and
according to the proposition \ref{singP} either the curve is
irreducible, or $x_1$ and $y_1$ are multiple zeros of the
discriminants $\tilde D_1(x)$ and $\tilde D_2(y)$, respectively, and
this singular point is unique. In both cases, we observe either the
case III of John's list of dynamical system behaviors and of a
breakdown of smoothness of the curve or a case of degeneration. In
this work we do not consider these cases. Some examinations were conducted in the work \cite{Lavr}.

Note also that the statement \ref{StatReduce} and proposition
\ref{StatRoots} of the subsection \ref{CanForms} imply

\begin{pr}\label{JohnPonc3}
For any two conics that are in one of generic dispositions 1)-5),
there exist real linear-fractional changes $R_1,R_2$ of real
parameters $x=R_1(\bar x)$, $y=R_2(\bar y)$, such that the
corresponding reduced biquadratic curve $C$ will have one of
the canonical form (i), (ii) or (iii) of the subsection
\ref{CanForms}.  Here the disposition 1) corresponds to the case
(iii), $-2<\hat b<2$
;\ the disposition 2) corresponds to the case (iii), $\hat b>2$; \
the disposition 3) corresponds to the case (iii), $\hat b<-2$;  \
the disposition 4) corresponds to the case (i), $\, c=1,\,\tilde
b>2$; \ the disposition 5) corresponds to the form (ii), $\,
c=-1$; the disposition 0) corresponds to the case (i), $a>0,\,
c=1$ with the condition of disappearence $\tilde b<2$ i.e.
$b^2<(a+1)^2$.
\end{pr}

\vskip 10pt

\subsection{Periodicity of Poncelet mappings.}\label{PerPonc}
Now we can apply the John mapping $T$ in order to construct two
sequences of points $x_n$ and $y_n$ on conics $A$ and $B$.
According to proposition \ref{JohnPonc1} of the subsection
\ref{PonJohn}, the Poncelet mapping $U_B$ gives us (and can be
obtained from) the mapping $U:\tilde C\to \tilde C$ that coincides
with the John mapping $T$ on the curve $C$ \re{2_2F}. Therefore
the periodic trajectories in the Poncelet problem correspond to
closure orbits of the John dynamical system. Thus, we obtain the main
result

\begin{pr}\label{PoncJohn}
The Poncelet problem in a generic setting is periodic if and only if
the John mapping on the corresponding biquadratic curve \re{Fdet} is
periodic, and then their periods are coincide.
\end{pr}

Explicit answers for the John mapping were given above, in propositions
\ref{StatRoots}, \ref{thetaJohn1}, \ref{etaJohn}. We obtain
\begin{pr}\label{PoncJohn1}
Dispositions 1)-5) of the previous subsection correspond to cases
1)-5) of the proposition \ref{StatRoots}. The conditions \re{per_theta}, 
\re{res1}-\re{res3} are criteria of periodicity of the Poncelet
mapping $U$ in corresponding cases.\end{pr}

Note that last proposition implies the statement of the big
Poncelet theorem.

Since the Poncelet problem is invariant under the arbitrary
projective transformation of the plane $(\xi, \eta)$ we can reduce
the conics $A$ and $B$ to some simpler shapes. We will consider the following possibilities:

(i) If we reduce $A$ and $B$ to concentric quadrics determined by
the equations
$$ \xi^2/a_1 + \eta^2/b_1 =1, \quad \xi^2/a_2 + \eta^2/b_2 =1.$$

Then the parameterization will take the form:
$$
\xi= 2a_1^{1/2} y/(1+y^2), \quad \eta= b_1^{1/2} (1-y^2)/(1+y^2)
$$
for the conic $A$ and
$$
\xi= 2a_2^{1/2} x/(1+x^2), \quad \eta= b_2^{1/2} (1-x^2)/(1+x^2)
$$
for the conic $B$. It can be easily verified that the polynomial
$Z(x,y)$ defined by \re{Fdet} can be reduced in this case to the
simplest Euler-Baxter form: $$ Z= x^2y^2 + 1 + a (x^2+y^2) + 2
b x y  
$$
with complex coefficients.

Hence, in this case, we have a simple form solution of the Poncelet problem, as above (see \re{M_n_Bax}): \be x_n=\sqrt{k}\:
\mbox{sn}(h(n+s_0);k), \quad y_n =\sqrt{k}\:
\mbox{sn}(h(n+s_0+1/2);k).
\ee

Parameters $a_i$ and $b_i$ can be specialized even further. For
instance, it is possible to choose $a_2=b_2=1$ reducing $B$ to a
unit circle. This choice corresponds to the so-called Bertrand model
of the Poncelet process \cite{Schoe}. If the second conic is an
ellipse $ \xi^2/a_1^2 + \eta^2/b_1^2 =1, \ a_1>1,\ b_1<1$, then the
formula \re{Fdet} gives a bounded curve of the form \be x^2 y^2
+ \frac{1+b_1}{1-b_1}(x^2+y^2) - \frac{4a_1}{1-b_1} xy + 1
=0.\lab{e_bax1}\ee In this case the points $x_n$ are isomorphic to
the godograph distribution of spins in the classical $XY$-chain
(see section \ref{spin} and \cite{GZ}). Another possible choice
$a_2-a_1=b_2-b_1$ corresponds to the confocal quadrics. In this
case the Poncelet problem is equivalent to the elastic billiard
\cite{Ves}, \cite {Laz}.

(ii) Let us consider a possibility of reducing the conics $B$ and
$A$ to two parabolas in the euclidean plane $(\xi, \eta)$. One of
the conics, say $B$, can be fixed by the chooing $\xi=x, \eta=x^2$,
whereas the other parabola remains arbitrary: $\xi=F_1(y)/F_0(y), \;
\eta=F_2(y)/F_0(y),$ where $F_0(y)=(ay+b)^2$ is square of a
linear function (this is a characteristic property of any
parabola). Then we obtain \be Z(x,y)=x^2F_0(y) -2x F_1(y) + F_2(y).
\lab{WP_Z}\ee Performing additional projective (complex)
transformation of the variable $y$, we can fix polynomials $F_1(y),
F_2(y),$ and $F_0(y)$ in such a way that the polynomial $Z(x,y)$
becomes symmetric in $x,y$. Then $Z(x,y)$ can be reduced to the
form \be Z=(xy + (x+y)y_0 + g_2/4)^2 - (x+y+y_0)(4xyy_0-g_3),
\lab{W_Z} \ee where $g_2,g_3$ are two remaining (arbitrary)
independent parameters of the polynomial $Z$.

In this case, we obtain the Poncelet points parameterized
by the Weierstrass function
$$
x_n = A_1 \wp(h(n+s_0)) + A_0 ,\quad y_n = A_1 \wp(h(n+s_0+1/2)) +
A_0.
$$
As we will see in subsection \ref{Toda}, the biquadratic curve of a similar form appears for the phase portrait of the elliptic solution of the
Toda chain.

Finally, we note that when $F_0=const$ and $F_1(y)$ is a linear
function in $y$, then the axis of the parabola $A$ is parallel to that of the parabola $B$. But then the curve
$Z(x,y)=0$ describes some arbitrary conics in the coordinates $x$ and
$y$ with the omitted $x^2y^2, x^2y,$ and $xy^2$ terms.

\subsection{Cayley determinant criterion}\label{Cayley} 
Let $A,B$ be some arbitrary conics, as described in previous section.
Recall that all tangents pass through points of the conic
$A$, whereas all vertices lie on conic $B$ (it is possible to
assume that the conic $A$ is located inside the conic $B$). Let $M_A$ and
$M_B$ be $3\times 3$-matrices describing conics (i.e.
corresponding quadratic forms) in the projective coordinates
$x_0,x_1,x_2$. This means that if the conic $A$ is a unit circle $x^2 +
y^2=1$, and the conic $B$ is a concentirc circle of radius $R$, then the
quadratic forms for $A,B$ take the form \be x_1^2 + x_2^2 -x_0^2 \quad
\mbox{and} \quad  x_1^2 + x_2^2 -R^2 x_0^2. \lab{2_circles} \ee
The corresponding matrices $M_A,M_B$ are diagonal:
$$M_A=diag(1,1,-1),  M_D=diag(1,1,-R^2).$$ Compute the
characteristic determinant \be F(z) = \det(A - z B). \lab{char_F}
\ee Clearly, $F(z)$ is a cubic polynomial. This polynomial is
invariant with respect to any similarity transformation $A \to
S^{-1} M_A S, \; B \to S^{-1} M_B S$ with a nondegenerate matrix
$S$. As is well known, for a pair of quadratic forms, in general,
apart from degenerate cases, there exist a transformation, which reduces
both forms to a diagonal form at the same time. The roots $z_i,
i=1,2,3$ of the polynomial $F(z)$ have a simple meaning. If the
matrix $M_B$ is reduced to identity matrix (i.e. $M_B =
diag(1,1,1)$), then $z_1,z_2,z_3$ are diagonal elements
(eigenvalues) of the matrix $A$.

The first step in the Cayley criterion is the Taylor expansion of the
square root of the polynomial $F(z)$: \be \sqrt{F(z)} = c_0 + c_1
z + c_2 z + \dots + c_n z^n + \dots \lab{Taylor_F} \ee

Then we compute the Hankel-type determinants:

\be H^{(1)}_p = \left | \begin{array}{cccc} c_3 & c_{4} & \dots &
c_{p+1}\\ c_{4}& c_{5} & \dots & c_{p+2}\\ \dots & \dots & \dots & \dots\\
c_{p+1} & c_{p+2} & \dots & c_{2p-1} \end{array} \right |, \quad
 p=2,3,4,\dots \lab{H1} \ee

and

\be H^{(2)}_p = \left | \begin{array}{cccc} c_2 & c_{3} & \dots &
c_{p+1}\\ c_{3}& c_{4} & \dots & c_{p+2}\\ \dots & \dots & \dots & \dots\\
c_{p+1} & c_{p+2} & \dots & c_{2p} \end{array} \right |, \quad
p=1,2,3,\dots \lab{H2} \ee

Then the Cayley criterion \cite{Berger}, \cite{GH} becomes:
\begin{st}
The trajectory of the Poncelet problem is periodic with the period equal to
$N$ if and only if $H^{(1)}_p=0$ for $N=2p$ and $H^{(2)}_p=0$ for
$N=2p+1$. For a modern proof of the Cayley criterion, see \cite{GH}.
\end{st}

{\bf Illustration.} Let us take once again the simplest case of two circles
\re{2_circles} with radii 1 and $R$. Then obviously, \be F(z) =
(z-1)^2(zR^2-1). \lab{F_circle} \ee The first nontrivial Taylor
coefficients of $\sqrt{F(z)}$ are $c_2 = R^2(R^2-4)/8$ and $c_3 =
R^4(R^2-2)/16$.  The case $c_2=0$ means that $R=2$, which corresponds to a
perfect triangle-shaped trajectory. The case $c_3=0$ means that $R=2^{1/2}$, which corresponds to a square-shaped trajectory.

\section{The Pell-Abel equation} \label{PAbel}
\subsection{Historical notes on the Pell-Abel equation and some
equivalent problems of calculus.} \label{PAbelH}

The Pell equation \be P^2-RQ^2=L\lab{Pell}\ee is the well-known
Diophantus equation, where for a given integer number $R$, which
is not a square, one seeks integers $P,Q,L$ satisfying the equation.

It was established that L. Euler had attributed this equation to Pell by misunderstanding. Actually the equation \re{Pell} arose and was studied in works of an Indian mathematician Brahmagupta as far back as 1000 years ago, i.e. before Euler (about 600 years earlier), and also in works of Euler's predecessors (P. Fermat was among them). Moreover, there is even a famous
problem by Archimedes related to a so-called Pell equation (Archimedes' cattle problem, see \cite{Archim}).

Theory of the Pell
equation is well-known (\cite{Mordell}); note only that the
standard theory of the equation \re{Pell} establishes a
connection with a continued fraction of the number $\sqrt{R}$.

The equation \re{Pell} in a ring of polynomials ${\R}[t]$ or
${\C}[t]$ of one variable with a constant $L$ is called the
Pell-Abel equation (it also appears under the names of "Abel equation" and "Pell equation for polynomials"). This equations arose in Abel's
work of 1826 year, in which he studied representation of the primitive
$\int{\rho(t)}/{\sqrt{R(t)}}\,dt$ through elementary functions; here
$\rho, R\ $ are polynomials
. N.H. Abel proved that if this primitive can be represented by
logarithm and rational functions of $t$ and $\sqrt{R}$, then one
can find polynomials $P,Q$ and a number $A$ such that \be
\int\frac{\rho}{\sqrt{R}}dt=A\ln\frac{P+\sqrt{R}Q}{P-\sqrt{R}Q}.\lab{AbIntG}\ee
Here the degree of $R\ $ is even: $\deg R=2m,$ $\deg\rho=m-1$ and
$\rho /A= 2P'/Q$. The main point of Abel's considerations was that
the polynomials $P$ and $Q$ satisfy the equation \re{Pell} with
$L=1$. And conversely, if the polynomials $P$ and $Q$ satisfy the
equation \re{Pell} with $L=1$, then the equality \re{AbIntG} holds
with $\rho = 2P'/Q$ and $A=1$.

Thus, the solvability of the Pell-Abel equation plays a role an
integrability criterion of the Abel differential. It is well-known
that Liouville and successors later showed that if the integral on the
left-hand-side part of \re{AbIntG} with some $\rho$ and $R$ can be
expressed via elementary functions, then the right-hand-side part must have
appearance \re{AbIntG} with some $P$ and $Q$ and the same $\rho$.
Note also that Abel gave one more criterion for the representation
\re{AbIntG}, namely, it was stated that the formula \re{AbIntG} holds,
iff a polynomial
continued fraction of the function $\sqrt{R}$ is periodic. Thus,
the solvability of the Pell-Abel equation also plays a role of a
criterion of periodicity of the continued fraction.

Let us consider one more classical problem -- Chebyshev's problem of a
search for a polynomial of least deviation. Let us consider a system
of $l$ closed intervals on the real axis
$I=[-1,1]\setminus\bigcup\limits_{j=1}^{l-1}(a_j,b_j)\ $ and the
following problem. One seeks a polynomial of the given degree
$n$ with a unit leading coefficient that has a least deviation on
the set $I$, i.e. one seeks a minimum of the functional $\|
t^n-P_{n-1}(t) \|_{C(I)}$. A general polynomial $P_{n-1}(t)$ runs
a finite-dimensional lineal space, and, hence, the problem can be
interpreted as a problem of search for an element in a
finite-dimensional lineal subspace of the Banach space, which is
nearest to the given element. But because the space ${C(I)}\ $ is not
reflexive, such nearest element does not necessarily exist.

Using a method that goes back to P. L. Chebyshev, existence of such
minimal polynomial was proven in \cite{AkhiezerLektsii}.

Further, if the polynomial $P_{n-1} $ is minimal on the set $I$,
then, possibly, it will also be minimal on some greater closed
subset $E\subset [-1,1]$.
Such a wider subset $E=[-1,1]$
$\setminus\bigcup\limits_{j=1}^{m-1}(\alpha_j,\beta_j)$ is called a
$n$-right extension of the set $I$ \cite{SodinYud}. If now one
takes a polynomial $R$ in the form
$R=(t^2-1)\prod\limits_{j=1}^{m-1} (t-\alpha_j)(t-\beta_j)$, then, as it was already proven, the solvability of the Pell-Abel equation
\re{Pell} with unknowns $P(t),Q(t),\,L=\mbox{const}>0$ is
equivalent to the statement, that the set $E$ is $n$-right extension of the set $I$. In addition, the polynomial $P$ gives a solution of the extremal
problem, and the number $\sqrt{L}$ is a least deviation (minimum
of the norm) \cite{SodinYud}.

It is interesting to note that, in this case, the subset $E$ is continuous
spectrum (which is absolutely continuous and two-valued) of some
bi-infinite Jacobi (3-diagonal) self-ajoint real periodic matrix in
the space $l^2$ if and only if the set $E$ is $n$-right or, equivalently, if the Pell-Abel equation \re{Pell} is solvable (see
references in \cite{SodinYud}). Note that the Chebyshev
problem will be obtained for the case of one interval $E=[-1,1],\
R=(1-t^2)$, and the Akhiezer's problem -- for the case of two
intervals ($\deg R=4$) and a polynomial of the fourth order
$E=[-1,a]\cup [b,1],a<b,\ R=(1-t^2)(t-a)(t-b)$. Note that there
exist several solvability criteria for the Pell-Abel equation with
a polynomial corresponding to the Akhiezer problem; among them
is the well-known Zolotarev porcupine (see, for instance,
\cite{Malyshev}).

\subsection{Connection between the Poncelet problem and the Pell-Abel equation}\label{PAPoncelet}
The solvability of the Pell-Abel equation \be
P(t)^2-R(t)Q(t)^2=L^2\lab{PellAbel}\ee with  the polynomial $R$ of an
even order, as was noted above, has several equivalent
formulations (\cite{SodinYud}). In the work of V.A. Malyshev
\cite{Malyshev}, there is a new solvability criterion given in an
algebraic form that we will formulate for a case, interesting for us,
of the fourth order $R=t^{4}+d_1t^{3}+...+d_{4}$. Let us expand
the root $\sqrt{R}$ into the Laurent series in some neighborhood of
infinity:
\be\sqrt{R}=\sum\limits_{j=-m}^{\infty}C_jt^{-j}\lab{sqrtR}\ee and
construct a determinant of the Hankel type: \be \Gamma_k = \left |
\begin{array}{cccc} C_1 & C_{2} & \dots &
C_{k}\\ C_{2}& C_{3} & \dots & C_{k+1}\\ \dots & \dots & \dots & \dots\\
C_{k} & C_{k+1} & \dots & C_{2k-1} \end{array} \right |, \quad
 k=1,2,3,\dots \lab{M1} \ee

{\it The Malyshev criterion} (\cite{Malyshev}) states:

\begin{st}
The Pell-Abel equation \re{Pell} with the polynomial $R$ of the fourth
order has, as a solution, some polynomials $P$ and $Q$ of orders $k+2$ and $k$ if and only if $\Gamma_k=0$.
\end{st}

Our observation is following. Let consider the Pell-Abel equation
\re{PellAbel} with $R$ of the fourth order and let $\lambda_1$ be
one of roots of the polynomial $R(t)$. Let us perform a shift of the
parameter $t\to t+\lambda_1$ (i.e. $t=\tilde t+\lambda_1$). A new Pell-Abel
equation \re{PellAbel} will also be solvable. Now let us apply the
Malyshev criterion to the new equation and let the coefficients
$C_j$ be coefficients of the Laurent series of the root
$\sqrt{R(t+\lambda_1)}$. Then, the equality $\Gamma_k=0$ is
necessary and sufficient for solvability of the equation
\re{PellAbel} with the given orders of polynomials $P,Q$. On the other
hand, let us perform the change of variable
$s=1/(t-\lambda_1),\ t=\lambda_1+1/s$ in the initial setting. Then we obtain \be
R(t)=F(s)/s^4\lab{RF}\ee with a polynomial $F(s)$ of the third
order and therefore $\sqrt{R}=\sqrt{F(s)}/s^2.$ Note that the
polynomial $R$ can be restored from $F$ by the inverse transformation:
\be R=t^4F(1/t)=t(z_1t-1)(z_2t-1)(z_3t-1),\label{PolynR}\ee where
$z_1,z_2,z_3-$ roots of the polynomial $F$ of the third order. The
polynomial $F$ can be generated by formula $F=\det (A-zB)$ and
matrixes $A=diag(z_1,z_2,z_3)$, $B=diag(1,1,-1)$ that are built by
quadratic forms $ z_1x_1^2 + z_2x_2^2 +z_3 x_0^2 \quad \mbox{and}
\quad x_1^2 + x_2^2 -x_0^2 $ of the conics \be z_1x^2 + z_2y^2 +z_3
=0,\ \  x^2 + y^2 =1.\label{conics}\ee

The expansion \re{Taylor_F} of the root $\sqrt{F(s)}$ gives us the
expansion $ \sqrt{F(s)}/s^2 = c_0/s^2 + c_1/s + c_2 + \dots + c_n
s^{n-2} + \dots $ which, after an inverse change, can be written as
$$ \sqrt{R(t+\lambda_1)} = c_0t^2 + c_1t + c_2 + c_3/t+ \dots +
c_n/t^{n-2} + \dots .$$ We see that $C_1=c_3,C_2=c_4,...$ and that
the determinant \re{H1} coincides with the determinant \re{M1},
i.e. the Caley criterion for the case of an even period with the
polynomial $F$ coincides with the Malyshev's criterion for the
equation \re{PellAbel} with the polynomial $R$. We have proven

\begin{pr} \lab{t_PA_Poncelet}
Let the Poncelet problem for the given conics $A$ and $B$ generate a
polynomial $F$ of the third order by formula \re{char_F} and a
polynomial $R$ by formula \re{RF}. Then the equation \re{PellAbel}
with polynomial $R$ of
the fourth order is solvable, if 
the Poncelet problem is periodic with an even period. Conversely,
for a given real polynomial $R$ of the fourth order without a constant
term, one can build a pair of conics $A$ and $B$, so that the
Poncelet problem (generated by polynomial $F$ by means of formulas
\re{char_F} and \re{RF}) is periodic with an even period, if the
equation \re{PellAbel} with polynomial $R$ is solvable.
\end{pr} 

Now the proposition \ref{etaJohn} gives us solvability criteria
for the Pell-Abel equation \re{PellAbel}.

\section{Criterions of uniqueness breakdown for the Dirichlet
problem and Ritt's problem}\label{breakdown} 

In this section we will give criteria of solution uniqueness for the
Dirichlet problem \re{SE},\re{NDP}
and show how to construct a 
system of some special solutions in the case of the non-unique the
Dirichlet problem for the string equation. We will describe a method
that can be applied to a slightly more general class of curves
than biquadratics.

We will consider the Dirichlet problem \re{SE},\re{NDP} in the
classical setting \re{UHDP}, if the curve $C$ is bounded, and in the
modified setting \re{MHDP}, if the curve $C$ is unbounded, because
we will need an application of the proposition \ref{unique}.

To remind, in subsection \ref{Jbc} 
we considered an even periodic meromorphic function $\phi(z)$ 
with some (complex) period $T$ and a complex curve $\tilde
C\subset \Bbb C^2$ described parametrically as \be x(t) =
\phi(t),\quad y(t) = \phi(t+\ve),\ \
t\in\Bbb C 
\ee with the property \re{CDom}, $\ve$ is a nonzero complex
parameter. We considered also a real curve $C=\tilde C\cap \Bbb
R^2$ given by means of a contraction of the functions $x,y$ on a
set
$S\in\Bbb C$. 
We saw that the curve $\tilde C$ is symmetric with respect to the
line $y=x$, the complex John mapping $T$ is equivalent to a shift
of the parameter $t$ by the step $-2 \ve$, and the periodicity
condition of the complex John mapping has the form: $ 2n\epsilon
=mT \lab{PerConOne}$ with some integer $n,m$.

We assume also (this is a very strong restriction, as we will see soon) that the function $\phi(z)$ possesses the nontrivial
multiplication property: \be \phi(n z) = R_n(\phi(z)), \;
n=1,2,3,\dots \lab{mult_prop} \ee where $R_n(z)$ is a family of
rational functions of the argument $z$ (by definition, $R_1(z)=z$,
but for other $n=2,3,\dots$ the expression and even the degree of
$R_n(z)$ can be non-obvious). We will consider the condition \be
\forall x\in\Bbb R,\ \ R_n(x)\neq \infty \label{Rnoninfty}\ee

Consider the complex setting \re{ComlSet} of the homogeneous Dirichlet
problem. Now the sufficient condition of the uniqueness of the
Dirichlet problem from proposition \ref{uniqueCom} requires the
transitivity of $T$ that implies that for any integers $n$ and
$m$, we have $2 n \ve \ne m T$.

Now let us assume that this condition is not fulfilled, i.e. for some
integers $N,M$, we have the condition \be 2 N \ve = M T
\lab{non_uniq} \ee

We will show that in this case the problem is, indeed, non-unique and we
will construct explicitly a system of explicit solutions \be
\Phi_n(x,y)=f_n(x)+g_n(y), \; n=1,2,\dots \lab{gen_string} \ee for
the string equation in the domain bounded by the curve $C$.

As the rational functions $R_n(z)$ are assumed to be non zero,
bounded on $\Bbb R^2$ and non-constant, then the function
$\Phi(x,y) = f_n(x)+g_n(y)= R_{2nN}(x) - R_{2nN}(y)$ is obviously
nonzero and smooth in $\Bbb R^2$. So we should verify the
Dirichlet boundary condition $\Phi(x,y)\equiv 0$ in all points of
the curve $C$. Indeed, for any point $t$ on the curve, we have
$$
f_n(x(t)) + g_n(y(t))= R_{2nN}(x(t)) - R_{2nN}(y(t)) = \phi(2nNt)
- \phi(2nNt - 2nN \ve) = 0,
$$
where we utilized properties \re{xy_phi}, \re{mult_prop} and
\re{non_uniq}. Hence, the function $\Phi_n(x,y)$ is identically
zero in all points on the curve $C$.  

We proved the following

\begin{pr} \lab{t_sols}
Under the condition \re{non_uniq}, for the curve $C$ we can choose
\be f_n(x) = R_{2nN}(x),\ \ \ g_n(y)= - R_{2nN}(y), \lab{fg_expl}
\ee where $R_n(z)$ are rational functions defined by
\re{mult_prop} with the property \re{Rnoninfty}. This will provide a
non-zero solution \re{gen_string} of the Dirichlet problem in the
complex setting \re{ComlSet}. If, in addition, the condition
\re{Rnoninfty} is fulfilled then the obtained solution satisfies the
modificed setting \re{MHDP} and, furthermore, if the curve is
bounded, then we will have a classical solution in the interior of the curve $C$, i.e. in the setting \re{UHDP}.
\end{pr}

First, we will illustrate how this theorem works in simplest case, when
the curve $C$ is an ellipse. In this case, we can choose the
parameterization \be x(t) = \cos(t), \quad y(t) = \cos(t+\ve)
\lab{par_ellips} \ee with some parameter $\ve$. Note that an
ellipse described by \re{par_ellips} has semiaxis $a=\sqrt{2}
\cos(\ve/2), b=\sqrt{2} \sin(\ve/2)$, one of which is inclined by
the angle $\pi/4$ with respect to axis $x$.

Now the John criterion \re{non_uniq} of solution uniqueness becomes
\be \ve \ne N \pi/M. \lab{J_ellips} \ee Otherwise, assume that
$\ve = N \pi/M$ for some positive integers $N,M$. The function
$\cos z$ satisfies all the needed conditions; moreover, we have \be
\cos(zn) = T_n(\cos(z)), \; n=1,2,\dots \lab{ChebT} \ee where
$T_n(x)$ are the Chebyshev polynomials of the first kind, hence
the proposition \ref{t_sols} gives us explicit solutions in this
case \be f_n(x)=T_{2Nn}(x), \; g_n(y)=-T_{2Nn}(y), \; n=1,2,\dots
\lab{Cheb_sol} \ee

Thus, in the case of nonuniqueness of the Dirichlet problem for an
ellipse we obtain a denumerable family of polynomial nonzero
solutions inside the ellipse.

Let $C$ be the biquadratic curve. To begin with, let us consider the
case \re{par13} of subsection \ref{Jbc}: $a>0,\ c=-1$. Here the
curve is parameterized via
$$
x=\phi(t)=\sqrt{\frac k {k'} } \: {cn}(t;k), \quad
y=\phi(t+\eta)=\sqrt{\frac k {k'}  } \: {cn}(t+\eta;k).
$$
An easy analysis shows that the last curve will be real for real $t$
and real $\eta$. For the elliptic cosine, there exist well-known
multiplication formulas
$$
cn(2mz)=T^1_m(cn^2\,z),\quad cn((2m+1)z)=cn\,z\,T^2_m(cn^2\,z)
$$
with some rational functions $T^1_n(z), T^2_n(z).$ Poles of the
function $cn(z)$ lie on the line  Im$\,t=K'$; therefore, the
condition \re{Rnoninfty} holds, the curve $C$ is bounded in this
case and the proposition \ref{t_sols} states that the condition
\re{res1} is sufficient for uniqueness breakdown for the Dirichlet
problem in the classical setting. Necessity of the condition \re{res1}
follows from the proposition \ref{unique} (for details, see
\cite{Telitsina}).

In a similar way, all other cases can also be considered (it was done in
works \cite{BurZhedan1},\cite{Bel },\cite{Telitsina}). Thus, we will
formulate the final result as follows:

\begin{pr}\label{CondBreak} The conditions
\re{res1}-\re{res3} are criteria of uniqueness breakdown for the
Dirichlet problem \re{SE},\re{NDP} in the corresponding cases
\re{par13}-\re{par a2} of the biquadratic curve in a canonical
form. The solutions are understood in the modified setting \re{MHDP}
of the Dirichlet problem, and for cases of the bounded curve $C\,-$
in the classical setting \re{UHDP}. Under the condition, there is a
denumerable set of linearly independed analytical solutions of the
homogeneous Dirichlet problem.
\end{pr}

Now we can consider the problem: what is the most general class of
even periodic functions $\phi(z)$ with the property
\re{mult_prop}? This problem has solved by Ritt in 1922
\cite{Ritt1}. The list of such possible functions includes two
generic cases:

(i) $\phi(z) = \frac{a \cos(z) + b}{c \cos(z) + d}$ with arbitrary
$a,b,c,d$;

(ii) $\phi(z)$ is a generic even elliptic function of the second degree,
i.e. $\phi(z) = \frac{a \wp(z) + b}{c \wp(z) + d}$, or,
equivalently, $\phi(z) = \kappa \: \frac{\sigma(t-e_1)
\sigma(t+e_1)}{\sigma(t-e_2) \sigma(t+e_2)}$ with arbitrary
parameters $a,b,c,d,\kappa, e_1,e_2$. Parameters $g_2,g_3$ of the
Weierstrass functions are arbitrary (equivalently, periods $2
\omega_1, 2\omega_2$ are arbitrary).

There are also three exceptional cases connected with elliptic
functions with certain restrictions imposed upon their parameters:

(iii) $\phi(z) = \frac{a \wp^2(z) + b}{c \wp^2(z) + d}$ with
arbitrary $a,b,c,d$ but with the restrictions that $g_3=0$ and $g_2$
is an arbitrary real. In this case, the fundamental parallelogram of
periods $2\omega_2, 2\omega_3$ is a simple square. Such case
corresponds to the so-called  {\it lemniscatic} elliptic functions.

(iv) $\phi(z) = \frac{a \wp^3(z) + b}{c \wp^3(z) + d}$ with
arbitrary $a,b,c,d$ but with restriction $g_2=0$ and $g_3$
arbitrary real. In this case, the period lattice has the hexagonal
symmetry. Such case corresponds to the so-called {\it equianharmonic}
elliptic functions.

For the equianharmonic case, Ritt found one more case when
$\phi(z)$ is a linear rational combination of the Weierstrass
function $\wp'(z)$, but in this case the function $\phi(z)$ is not
even, and we can omit this case.

The case (ii) corresponds to our problem in biquadratic. Thus, in
the case of non-uniqueness, we have the solution \be f_n(x)=
R_{2nN}(x), \quad g_n(y)=-R_{2nN}(y), \lab{sol_ell} \ee where
$R_n(z)$ are rational functions defined as \re{mult_prop}. Note
that in contrast to the trigonometric case (i.e.
$\phi(z)=\cos(z)$), there exist no simple explicit expressions for
$R_n(z)$ for arbitrary $n$.

What about cases (iii) and (iv)? In these cases, the curve $C$
is described by an algebraic equation of degree greater
than four. Consider, for example, the case (iii) and assume for simplicity
that $\phi(z)=\wp^2(z)$. We seek an algebraic equation
$f(x,y)=0$ for the curve $C$ described parametrically as
$x(t)=\wp^2(t), \; y(t)=\wp^2(t+\epsilon)$. This equation can be
easily derived from the addition theorem for the Weierstrass function.
We will not write it down explicitly, but, instead, we will only mention that the degree of this equation is eight. Similarly, for the case (iv), we
obtain the algebraic curve described by a polynomial $F(x,y)$ of
the overall degree of 12. But we should note that the John condition of
intersection at no more than two points with vertical and
horizontal lines is not satisfed and here we do not have a
sufficient condition of solution uniqueness for boundary value
problems.

\section{Related problems of mathematical physics}

\subsection{Classical Heisenberg XY spin chain}\label{spin} 
There is an interesting relation to the classical Heisenberg
$XY$ spin chain \cite{GZ}, which is a system of two-dimensional unit
vectors ("spins") $\vec r_n=(x_n;y_n), \: |\vec r_n|=1$ with the
energy of interaction given as \be E=\sum_{n=0}^{N-1} (\vec r_n, J \vec
r_{n+1}), \lab{EXY} \ee where $J =diag(J_1,J_2)$ is a 2-by-2 diagonal
matrix. The problem here is to find static solutions that
provide a local extremum to the energy $E$. As was shown in
\cite{GZ}, this problem is equivalent to finding solutions of the
systems of non-linear vector equations in the form: \be (\vec r_n, J(\vec
r_{n-1}+ \vec r_{n+1})) =0, \quad n=1,2,\dots,N-1. \lab{eq_XY} \ee

Below we will select cases of a closed chain, in which $\vec
r_0=\vec r_N$ for some $N$. Based on this, it will be assumed that
$J_1=1,J_2=j>1$. Among all the solutions of \re{eq_XY}, we choose the
so-called regular solutions \cite{GZ}, satisfying the condition: \be \vec
r_{n-1}+ \vec r_{n+1} \ne 0. \lab{reg_XY} \ee Then it is possible
to show that the scalar product \be W=(\vec r_n, J\vec r_{n+1})
\lab{W_XY} \ee does not depend on $n$ and, hence, it can be considered
as an integral of the system \re{eq_XY}. It is sufficient to construct
general regular solutions \cite{GZ}, which will be of two types. The choice
of the solutions depends on the value of the integral $W$. If $|W|<1/j$ then
\be x_n =cn(q(n-\theta);k),\quad y_n =sn(q(n-\theta);k),
\lab{sol1_XY} \ee where parameters $k,q$ can be found from \be
dn(q;k)=1/j, \quad k^2=\frac{1-j^{-2}}{1-W^2}. \lab{kq_1} \ee

If $1/j<|W|<1$ then \be x_n =dn(q(n-\theta);k),\quad y_n =k\:
sn(q(n-\theta);k), \lab{sol2_XY} \ee where \be cn(q;k)=1/j, \quad
k^2=\frac{1-W^2}{1-j^{-2}}. \lab{kq_2} \ee

In both cases, the parameter $\theta$ is an arbitrary real number
depending on an initial condition. If the chain is periodic, then we will
have \be qN=4Km_1 + 2iK'm_2 \lab{per_XY} \ee which coincides with
\re{per_eta}. The reason for such coincidence can be described as follows.

Let us consider an equation of the following integral \be x_n x_{n+1} + j^{-1} y_n
y_{n+1} =W \lab{W_XY1} \ee with a fixed value $W$. The equation
\re{W_XY1} can be reduced to the algebraic form \re{E_Bax} by a
standard substitution (stereographic projection from a unit
circle to a line):
$$
x_n = \frac{1-u_n^2}{1+u_n^2}, \quad y_n = \frac{2u_n}{1+u_n^2}$$
It is easy to see that the variables $u_n, u_{n+1}$ lie on the
Euler-Baxter biquadratic curve
$$
u_n^2 u_{n+1}^2 + 1 + a (u_n^2+ u_{n+1}^2) + b u_n u_{n+1}=0
$$
with parameters $a,b$ simply related to the "physical" parameters
$j,W$.

Then it can be easily verified that finding solutions (step-by-step)
of equations \re{eq_XY} for the regular solutions is equivalent
to finding points $M_2,M_3,\dots M_N$ for the John algorithm. Note
that arguments similar to those in the formula \re{per_theta} ) allow concluding that the parameters $q$ and $\theta$ must be real; hence, the
condition \re{per_XY} takes the form \be {q\over4K}={m_1\over N}.
\lab{per_XY2} \ee

Thus, static regular solutions for the closed finite classical
Heisenberg $XY$-chain are equivalent to periodic solutions of the
John algorithm for the Euler-Baxter biquadratic curve, or,
equivalently, to finding periodicity condition of the Poncelet
process for a unit circle $x^2+y^2=1$ and a concentric conic
$x^2/\xi_1 + y^2/\xi_2=1$. This choice of conics corresponds to
the Bertrand model of the Poncelet process \cite{Schoe}. This
means that there exists equivalence between static periodic
solutions of the Heisenberg $XY$-chain and the Bertrand model of
the Poncelet process.



\subsection{Elliptic solutions of the Toda chain and the biquadratic
curve} \label{Toda}
Now we consider the Toda chain that is a discrete dynamical system
consisting of two sets $u_n(t), b_n(t)$ of complex variables
depending on continuous parameter $t$ and discrete parameter $n=0,
\pm 1, \pm 2, \dots$. The equations of motion are \cite{Toda}
\be \dot b_n = u_{n+1}-u_{n}, \quad \dot u_n = u_n (b_{n} -
b_{n-1}) \lab{Toda_eq} \ee The Toda chain is one of the simplest and most
and famous completely integrable discrete dynamical systems
(for the history of this model and a review of different approaches, see
\cite{Toda}). Among numerous explicit solutions, there exist
so-called "elliptic waves" constructed for the first time by Toda himself
\cite{Toda}. We will give the Toda elliptic solutions in somewhat
different form which is more convenient for applications (Toda
used Jacobi elliptic functions, whereas we utilize the Weierstrass
functions).

\begin{pr}
An elliptic solution for the unrestricted Toda chain can be presented
in the from \ba &&b_n = \omega \Bigl(\zeta (\omega t - p(n+1) + r)
-  \zeta (\omega t -
p n +r)   \Bigr) - \lambda, \nonumber \\
&&u_n = \omega^2 \Bigl(\wp(p) - \wp(\omega t - p n + r) \Bigr) =
\nonumber \\
&&\omega^2 \: \frac{\sigma(\omega t - p(n+1) +r ) \, \sigma(\omega
t - p(n-1) +r)}{\sigma^2(p)\,  \sigma^2(\omega t - p n +r)}
\lab{w_sol} \ea with arbitrary parameters $\omega,p,r,\lambda$ and
arbitrary invariants $g_2,g_3$. Here $\wp(z), \sigma(z),$
$\zeta(z)$ are standard Weierstrass functions defined as in
\cite{WW}.
\end{pr}

{\it Proof}. The Toda chain equations \re{Toda_eq} are verified
directly by substitution using well-known formulas \cite{WW}:
$$
\frac{d}{d z} \zeta(z) = - \wp(z), \quad \frac{d}{d z} \lg \,
\sigma(z) = \zeta(z), \quad \wp(x) - \wp(y) =
\frac{\sigma(y-x)\sigma(x+y)}{\sigma^2(x) \sigma^2(y)}
$$
Strictly speaking, the variable $b_n(t)$ is inessential - it can
be eliminated from the system \re{Toda_eq}. Thus, only the variable
$u_n(t)$ is considered as a "true" Toda chain variable (for
details, see \cite{Toda}).

Now we can construct the phase portrait for this variable. By the
phase portrait, we will assume the plot constructed on the plane
$x,y$, if one take points $P_0,P_1,P_2, \dots$ with the coordinates
$P_0=(u_0,u_1), \; P_1=(u_1,u_2), \dots, P_n=(u_n, u_{n+1})$. The
phase portrait is an indicator of integrability: if the system is
integrable, then the points $P_i, i=0,1,\dots$ will fill some smooth curve.
Otherwise, these points are distributed quasi-stochastically
(the so-called "stochastic web" \cite{Zasl}). In our case, the variable
$u_n(t)$ is given explicitly by \re{w_sol},
$u_{n+1}(t)=u_n(t-p/\omega)$ and, hence, the points
$P_n=(u_n,u_{n+1})$ fill the symmetric biquadratic curve in its
canonical form \re{WP_curve} by the proposition \ref{www}. The
complex John algorithm in this biquadratic curve is equivalent to
passing from the point $P_n$ to point $P_{n+1}$ and then to point
$P_{n+2}$. Note that the time parameter $t$ describes a smooth
motion along this curve (Hamiltonian dynamics), whereas the John
algorithm describes a "discrete motion". The periodic case is
defined by a positive integer $N$ such that $u_n(t) = u_{n+N}(t)$.
From \re{w_sol} it can be seen that the periodicity condition is
equivalent to \be pN= 2 m_1 \omega_1 + 2 m_2 \omega_2
\lab{per_Tod} \ee with some integers $m_1, m_2$. An important
property of the elliptic Toda solutions lies in the fact that the periodicity property \re{per_Tod} does not depend on the time parameter $t$, i.e. if a periodic condition takes place for one value of $t$, then
it also holds for all others values of $t$. In terms of the
(complex) John algorithm, that means that a period of a point
following the John mapping does not depend on a choice of this
initial point $P_0=(u_0(t), u_1(t))$ on the curve.

Comparing obtained formulas for $u_n(t)$ with \re{WP_Z} and
\re{W_Z}, it can noticed that the elliptic solutions of the Toda chain
give an interesting illustration of the Poncelet theorem: periodic
solutions \re{per_Tod} of the Toda chain correspond to periodicity
of the Poncelet process on two parabolas.

\subsection{Elliptic grids in the theory of rational interpolations}
In this section we will describe briefly an interesting connection of
the Poncelet problem with theory of the so-called admissible grids for
biorthogonal rational functions. This subject is important in the
theory of rational interpolations. For further details and
relations with the theory of special functions, see, e.g.
\cite{Magnus}, \cite{SZ1}, \cite{SZ_Ky}, \cite{SZ_Ram},
\cite{ZheR}, \cite{Zhe_PI}.

Let us consider a set of rational functions $R_n(x)$ of the
order $[n/n]$, which means that $R_n(x)$ are given by ratios of
two $n$-th degree polynomials in $x$. It will be assumed that all rational functions $R_n(x)$ have only simple poles.

We take $\alpha_1, \alpha_2, \dots, \alpha_n$ as $n$ distinct
prescribed positions of the poles of $R_n(x)$. Then $R_n(x)$ can
be written down through a sum of partial fractions \be R_n(x) = t_0^{(n)} +
\sum_{i=1}^n \frac{t_i^{(n)}}{x-\alpha_i} \lab{pol_ex} \ee with
the coefficients $t_i^{(n)}, \; i=1,2,\dots,n,$ playing the role
of residues of $R_n(x)$ at the poles $\alpha_i$. The coefficient
$t_0^{(n)}$ can be interpreted as $\lim_{x \to \infty} R_n(x)$.

Let $x(s)$ be a "grid", i.e. a function in the argument $s$. We
would like to construct the so-called lowering operator ${\mathcal
D}_{x(s)}$ in the space of rational functions $R_n(x)$ defined on
this grid in the following way. As a definition of the
lowering operator ${\mathcal D}_{x(s)}$, we take a divided difference
operator in the parameterizing variable $s$, which obeys the
following properties:

i) the grid $x(s)$ is a meromorphic function of $s\in\mathbb C$,
which is invertible in some domain of the complex plane;

(ii) for any function $F(s)$ one has
$$
{\mathcal D}_{x(s)} F(s) = \chi(s) (F(s+1)-F(s)),
$$
where $\chi(s)$ is some function to be determined;

(iii) ${\mathcal D}_{x(s)} R_1(x) = const$, where $R_1(x)$ is an
arbitrary rational function of the order $[1/1]$ with the only
pole at $x=\alpha_1$;

(iv) the operator ${\mathcal D}_{x(s)}$ transforms any rational
function $R_n(x)$ with the {\it prescribed} poles  $\alpha_1,
\dots, \alpha_n$ to a rational function  $\t R_{n-1}(y(s))$ of the
adjacent grid $y(s)$ with some other poles $\beta_1, \dots,
\beta_{n-1}$;

(v) the operator ${\mathcal D}_{x(s)}$ is ``transitive": for any
nonnegative integer $j$, the operator ${\mathcal D}_{x(s)}^{(j)}$
defined as
$$
{\mathcal D}_{x(s)}^{(j)} F(s) = \chi_j(s) (F(s+1)-F(s))
$$
with some function $\chi_j(s)$ transforming any rational function
$R_n(x)$ with the {\it prescribed} poles $\alpha_{j+1}, \dots,
\alpha_{j+n}$ to a rational function  $\t R_{n-1}(y(s))$ with the same adjacent grid $y(s)$ and the sequence of poles
$\beta_{j+1}, \dots, \beta_{j+n-1}$;

(vi) we assume that the poles are nondegenerate: there are
infinitely many distinct values of $\alpha_n$ and $\beta_n$ and
$\alpha_n \ne \alpha_{n+1},  \alpha_{n+2}$ and, similarly,
$\beta_n \ne \beta_{n+1},  \beta_{n+2}$ for all $n$.

An important restriction is the condition of independence of
${\mathcal D}_{x(s)}$ on the order $n$ of a rational function. The
problem is to deduce the functions $\chi_j(s)$ and $x(s), y(s)$
from the given set of requirements.

From the properties (i)-(iii) one can easily find
\begin{eqnarray}\nonumber
&& \chi(s)= \left(\frac{1}{x(s+1)-\alpha_1} -
\frac{1}{x(s)-\alpha_1}\right)^{-1}
\\&&\makebox[2em]{}
=\frac{(x(s)-\alpha_1)(x(s+1)-\alpha_1)}{x(s)-x(s+1)}.
\lab{chi}\end{eqnarray} The function $\chi(s)$ is defined up to an
inessential constant multiplier. Note that from (ii), we have
${\mathcal D}_{x(s)} R_0(z)=0$.

The most non-trivial problem in the construction of the operator
${\mathcal D}_{x(s)}$ consists in establishing the properties
(iv)-(v).

This problem was solved in \cite{SZ_Ram}. It appears that the grid
$x(s)$ as well as the grids $y(s),\alpha_s, \beta_s$ should belong
to the class of elliptic grids. This means that they should
satisfy the biquadratic equation \ba&&A_1 x^2(s+1) x^2(s) + A_2
x^2(s) x(s+1) + A_3 x^2(s+1) x(s) + A_4 x(s+1) x(s) + \nonumber
\\&&A_5 x^2(s+1) + A_6 x^2(s) + A_7 x(s+1) + A_8 x(s) + A_9 =0
\lab{bi_grid} \ea with some constants $A_i, i=1,2,\dots,9$. As we
already know, this equation can be parameterized in terms of the
elliptic functions, and we, thus, arrive at elliptic grids $x(s)$.
Thus, the elliptic grids $x(s), y(s)$ are the most general grids to
provide existence of the lowering operator for rational functions.

In the theory of the rational (also called sometimes as the Cauchy-Jacobi
or Pad\'e) interpolation, these grids appear naturally for some
class of self-similar solutions.

The Cauchy-Jacobi interpolation problem (CJIP) for the sequence
$Y_j$ of (nonzero) complex numbers can be formulated as follows
\cite{BGM}, \cite{Mein}. Given two nonnegative integers $n,m$,
we choose a system of (distinct) points $x_j, j=0,1, \dots, n+m$ on
the complex plane. We seek polynomials $Q_m(x;n), P_n(x;m)$
of degrees $m$ and $n$, respectively, such that \be Y_j=
\frac{Q_m(x_j;n)}{P_n(x_j;m)}, \quad j=0,1, \dots n+m \lab{CJP}
\ee (in our notation we stress that, for example, the polynomial $Q_m(x;n)$, while being of degree $m$ in $x$, depends on $n$ as a parameter).

It may be that a solution of the CJIP does not exist. In this
case, it is reasonable to consider a {\it modified} CJIP in the form: \be Y_j
P_n(x_j;m)-Q_m(x_j;n)=0, \; j=0,1,\dots,n+m, \lab{MCJ} \ee where the
polynomials $P_n(x;m), Q_m(x;n)$ can now be unrestricted. The
problem \re{MCJ} always has a nontrivial solution. In an exceptional
case, if the system \re{CJP} has no solutions, some zeroes of the
polynomials $P_n(z;m)$ and $Q_m(z;n)$ coincide with interpolated
points $x_s$. Such points, in this case, are called unattainable
\cite{Mein}.

The CJIP is called {\it normal}, if the polynomials $Q_m(x;n),
P_n(x;m)$ exist for all values of $m,n=0,1,\dots$, and polynomials
$Q_m(x;n), P_n(x;m)$ have no common zeroes. This means, in
particular, that the polynomials $Q_m(x;n), P_n(x;m)$ have no roots,
coinciding with interpolation points, i.e. \be Q_m(x_j;n) \ne 0,
\; P_n(x_j;m) \ne 0, \quad j=0,1,\dots n+m \lab{n_deg} \ee In a
special case, in which there exists an analytic function $f(z)$ of a
complex variable such that $f(x_j)=Y_j$, the corresponding CJIP is
called the multipoint Pad\'e approximation problem \cite{BGM}.

It is possible to show that the Cauchy-Jacobi interpolation
problem is equivalent to the theory of biorthogonal rational functions
\cite{Zhe_PI}. The elliptic solutions (for the first time obtained in work
\cite{SZ1}) of this problem appear naturally, if the interpolated
function $f(z)$ satisfies the so-called discrete Riccati equation
\cite{Magnus}. Geometric interpretation of obtained elliptic grids
and their connection with the Poncelet problem can be found in
\cite{Magnus} and \cite{SZ_Ram}.

Note that if terms of degree $>2$ are absent in \re{bi_grid} (i.e.
$A_1=A_2=A_3=0$), then a corresponding grid is degenerated to the
so-called Askey-Wilson grid \cite{AW}, which is the most general
grid for orthogonal polynomials, satisfying a linear second-order
difference equation \cite{ViZhe}. From geometrical point of view,
the Askey-Wilson grids correspond to the John algorithm for the
second-degree curves (i.e. ellipsis, hyperbola or parabola)
\cite{Magnus}.

\bb{99}

\bi{Akhiezer} N.I. Akhiezer, {\it Elements of the Theory of
Elliptic Functions}, 2nd edition, “Nauka”, Moscow, 1970.
Translations Math. Monographs {\bf 79}, AMS, Providence, 1990.

\bi{AkhiezerLektsii}  N.I. Akhiezer, {\it Lectures on
approximation theory}, Nauka, M., 1965 (In Russian).

\bi{Alex1} R.A. Alexandrjan, On the Dirichlet problem for the
string equation and on completeness of a system of function in a
disk.-- Doklady AN USSR. 1950, 73, No.5 (In Russian).

\bi{Alex2} R.A. Alexandrjan, Spectral properties of operators
generated by systems differential equations of Sobolev type, Trudy
Mosc. Math. Obshchestva  9(1960), pp.455-505. (In Russian).

\bi{Alex3} G.S. Akopyan, R.A. Aleksandryan, On the completeness of
a system of eigen- and vector-polynomials of a linear differential
operator pencil in ellipsoidal domains, Dokl. Akad. Nauk Arm. SSR,
V.86, No.4, pp. 147-152 (1988). (In Russian).

\bi{Arnold} V. I. Arnold, {\it Small demominators. I}, Izvestija
AN SSSR, serija matematicheskaja, 25(1961), {\bf 1}, pp.21-86.

\bi{AW} R.~Askey and J.~Wilson, {\it Some basic hypergeometric
orthogonal polynomials that generalize Jacobi polynomials}, Mem.
Amer. Math. Soc. {\bf 54}, No. 319, (1985), 1-55.

\bibitem{BGM} G.A. Baker, P. Graves-Morris, {\it Pad\'e approximants.
Parts I and II.}. Encyclopedia of Mathematics and its
Applications, {\bf 13, 14}. Addison-Wesley Publishing Co.,
Reading, Mass., 1981.

\bi{BatErd} H. Bateman and A. Erd\'elyi, {\it Higher
transcendental functions. 3}, McGraw-Hill, New York, 1955, Bateman
manuscript project.

\bi{Bax} R. Baxter, {\it Exactly Solvable Models in Statistical
Mechanics}, London, Academic Press, 1982.

\bi{Bel }  M. V. Beloglyadov, {\it On the Dirichlet problem for
the vibrating string equation in domain with a bi-quadratic
boundary}, Trudy IAMM NASU, V.14, 2007, pp. 14-29 (In Russian).

\bi{Ber}  Yu. M. Berezanskii, {\it Expansion by eigenfunctions of
selfadjoint operators}, Naukova Dumka, Kiev, 1965.(In Russian)

\bi{Berger} M. Berger, {\it G\'eom\'etrie}, CEDIC, Paris, 1978.

\bi{BD} D. Bourgin, R. Duffin, {\it The Dirlchlet problem for the
vibrating string equations}, Bull.Am.Math.Soc., 1939, v.45,
pp.851-858.

\bibitem{BurAlgB}  V. P. Burskii, {\it On solution uniqueness of some boundary value
problems for differential equations in domains with algebraic
boundary}, Ukr. math. journal {\bf 45} (1993), No.7, pp. 993-1003.

\bi{BurMomProblem}  V. P. Burskii, {\it On boundary value problems
for differential equations with constant coefficients in a plane
domain and a moment problem}, Ukr. math. journal, {\bf 48},
(1993), No.11, pp. 1659-1668.

\bi{BurBook}  V. P. Burskii, {\it Investigation msthods of
boundary value problems for general differential equations}, Kiev,
Naukova dumka, 2002 (In Russian).

\bi{BurZhedan1}  V. P. Burskii, A. S. Zhedanov {\it On Dirichlet
problem for string equation, Poncelet problem, Pell-Abel equation,
and some other related problems}, Ukr. math. journal, {\bf 58}
(2006), No. 4, pp. 487-504.

\bi{BurZhedan2} V.P. Burskii, A.S. Zhedanov, {\it Dirichlet and
Neumann problems for string equation, Poncelet problem and
Pell-Abel equation} Symmetry, Integrability and Geometry: Methods
and Applications, 2006, V. 2, rec.No: 041.

\bi {BurZhedan3} V. P. Burskii, A. S. Zhedanov, {\it Boundary
value problems for string equation, Poncelet problem, and
Pell-Abel equation: links and relations}, Contemporary
Mathematics. Fundamental Directions, 2006, {\bf 16}, pp. 5–9.

\bi{Zasl} A.A. Chernikov, R.Z. Sagdeev,  G.M. Zaslavsky, G. M.
{\it Stochastic webs. Progress in chaotic dynamics}. Phys.D 33
(1988), no. 1-3, 65–76.

\bibitem{EK} O. Egecioglu and C.K. Koc, {\it  A fast algorithm for rational
interpolation via orthogonal polynomials}, Math. Comp. {\bf 53}
(1989), pp. 249--264.

\bibitem{Bat} A. Erdelyi, W. Magnus, F. Oberhettinger, F.G.
Tricomi, {\it Higher Transcendental Functions. I}, McGraw-Hill,
New York, 1953 Bateman manuscript project.

\bibitem{Fokin} M.V. Fokin, Solvability of the Dirichlet problem for
the string equation, Doklady AN SSSR, 1983, V. 272, No. 3,
pp.801-805 (in Russian).

\bi{FrRag} J.P.Francoise and O.Ragnisco, {\it An iterative process
on quartics and integrable symplectic maps}, in "Symmetries and
integrability of difference equations", P.A.Clarkson and
F.W.Nijhoff eds., Cambridge University press, 1998.

\bi{GZ} Ya.I. Granovskii and A.S. Zhedanov, {\it Integrability of
the classical $XY$-chain}, Pis'ma to Zh. Exp. Theor. Phys. {\bf
44} (1986), pp. 237--239 (Russian).

\bi{GH1} P. Griffiths and J. Harris, {\it Poncelet theorem in
space}, Comment. Math. Helvetici, {\bf 52} (1977), pp. 145--160.

\bi{GH} P. Griffiths and J. Harris, {\it On a Cayley's explicit
solution to Poncelet's porism}, Enseign. Math. (2), {\bf 24}
(1978), pp. 31--40.

\bi{GHP} P. Griffiths and J. Harris, {\it Principles of algebraic
geometry}, v. I,II, John Wiley and Sons, Inc., 1978.

\bi{Had} J. Hadamard, {\it Equations aux derivees partielles},
L`Enseignment Mathematique, Vol. 36(1936), pp. 25-42.

\bi{Halphen} G.H.~Halphen. {\it Trait\'e des Fonctions Elliptiques
et de Leures Applications}, II, Gauthier–Villar, Paris (1886).

\bi{Hub} A. Huber, {\it Erste Randwertaufgabe fur geschlossene
Bereiche bei der Gleichung $U_{xy}=f(x,y)$}, Monatshefte f\"ur
Mathematik und Physik, {\bf 39} (1932), pp. 79-100.

\bi{Ince} E.L.~Ince, {\it Ordinary differential equations.}

\bi{John} F.John, {\it The Dirichlet problem for a hyperbolic
equation}, Am.J.Math. {\bf 63} (1941), pp. 141--154.

\bi{IatrouRob} A. Iatrou, J.A.G. Roberts, {\it Integrable mappings
of the plane preseving biquadratic invariants curves II},
Nonlinearity {\bf 15} (2002), pp. 459-489.

\bi{IatrouArx} A. Iatrou,{\it Real Jacobian Elliptic Function
Parameterization for a Genuinely Asymmetric Biquadratic Curve},
arXiv: nlin. SI/0306051 v1 25 Jun 2003.

\bi {Kerawala} S. M. Kerawala,  {\it Poncelet Porism in Two
Circles}. Bull. Calcutta Math. Soc. {\bf 39}, pp. 85-105, 1947.

\bi{King} J.L. King, {\it Three problems in search of a measure},
Amer. Math. Monthly {\bf 101} (1994), pp. 609-628.

\bi{Lavr} M.M.~Lavrent'ev, {\it Mathematical problems of
tomography and hyperbolic mappings}, Sib. Math. J., {\bf 42}
(2001), No.5, pp. 916--925.

\bi {Laz} V.F. Lazutkin, {\it KAM Theory and Semiclassical
Approximation to Eigenfunctions}, Springer Verlag, Berlin,
Hei-delberg, New York (1993), Ergebnisse der Mathematik und
ihrer Grenzgebiete: 3. Folge, Band 24. 

\bi{Magnus} A.Magnus, {\it Rational interpolation to solutions of
Riccati difference equations on elliptic lattices}. Preprint
http://www.math.ucl.ac.be/membres/magnus/

\bi{Malyshev} V.A. Malyshev, {\it Abel equation}, Algebra and
analysis, {\bf 13} (2001), No. 6, pp. 1-55. (In Russian)

\bibitem{Mein} J. Meinguet, {\it On the solubility of the Cauchy interpolation problem}.
1970 Approximation Theory (Proc. Sympos., Lancaster, 1969) pp.
137--163. Academic Press, London.

\bibitem{Mordell} L.J. Mordell, Diophantine equations, Academic Press, 1969.

\bibitem {Nit}
Z. Nitecki,{\it Differentiable dinamics, MIT Press, Cambridge Mass
- London, 1971}

\bi{Ovs} S.G. Ovsepjan, On ergodisity of continuous automorphizms
and solution uniqueness of the Dirichlet problem for the string
equation. II.-- Izv. AN Arm.SSR. 2(1967), No.3, pp.195-209.

\bibitem{Pta}  B. Yo. Ptashnik, {\it Incorrect boundary value problems for
differential equations with partual derivatives}, Kiev, Naukova
dumka, 1984.(In Russian)

\bi{Ritt1} J.F.Ritt, {\it Periodic functions with a multiplication
theorem}. Trans. Amer. Math. Soc. {\bf 23} (1922), no. 1, pp.
16--25. 30E99

\bi{Schoe} I.J.Schoenberg, {\it On Jacobi-Bertrand's proof of a
theorem of Poncelet}. Studies in pure mathematics, To the Memory
of Paul Turan, 623--627, Birkhдuser, Basel, 1983.

\bi{SodinYud} L.M. Sodin, P.M. Yuditskii, {\it Functions least
deviating from zero on closed sets of real axis},
 Algebra and analysis, {\bf 4}, No.2, с.1-61.(In Russian)

\bibitem{SZ1} V.Spiridonov and A.Zhedanov, {\it Spectral transformation chains
and some new biorthogonal rational functions}, Commun. Math. Phys.
{\bf 210} (2000), 49--83.

\bibitem{SZ_Ky} V.P.Spiridonov and A.S.Zhedanov, {\it To the theory of
biorthogonal rational functions}, RIMS Kokyuroku {\bf 1302}
(2003), 172--192.

\bi{SZ_Ram} V.Spiridonov and A.Zhedanov, {\it Elliptic grids,
rational functions, and Pad\'e interpolation}, Ramanujan J., {\bf
13}, No. 1--3 (2007), 285--310.

\bibitem{Stieltjes} T.Stieltjes, {\it Sur l'\'equation d'Euler},
Bul.Sci.Math., Paris, s\'er. 2, {\bf 12} (1888), 222--227.

\bi{Telitsina} A.A. Telitsyna {\it The Dirichlet problem for wave
equation in plane domain with biquadratic boundary,} Trudy IAMM
NASU, V.13, 2007, pp. 198-210 (In Russian).

\bibitem{Toda} M. Toda, {\it Theory of nonlinear lattices},
Springer Series in Solid-State Sciences, vol. {\bf 20},
Springer-Verlag, Berlin, 1989.

\bi{Ves} A.P. Veselov, {\it Integrable systems with discrete time
and difference operators}, Functional Analysis and its
Applications {\bf 22} (1988), 1--13 (Russian).

\bi{VesUs} A.P. Veselov, {\it Integrable maps}, Russian Math.
Surveys {\bf 46} (1991), 1--51.

\bi{WW} E.T. Whittacker, G.N. Watson, {\em A Course of Modern
Analysis}, Cambridge, 1927.

\bi{ViZhe}  L.Vinet and A.Zhedanov,  {\it Generalized Bochner
theorem: characterization of the Askey-Wilson polynomials},
J.Comput.Appl.Math. {\bf 211} (2008), 45 –- 56.

\bi{Zell} T.I. Zelenjak, Selected topics of quality theory of
equations with partial derivatives.-- Novosibirsk: NGU, 1970.
(In Russian)

\bibitem{ZheR} A. Zhedanov, {\it Biorthogonal rational functions and the
generalized eigenvalue problem}. J. Approx. Theory {\bf 101}
(1999), 303--329.

\bi{Zhe_PI} A.Zhedanov, {\it Pad\'e interpolation table and
biorthogonal rational functions}, Proceedings of the Workshop on
Elliptic Integrable Systems November 8-11, 2004, Kyoto, Rokko
Lectures in Mathematics, No. 18, 323--363.
http://www.math.kobe-u.ac.jp/publications/rlm18/20.pdf

\bibitem{Archim}
$http://en.wikipedia.org/wiki/Archimedes\%27\_ cattle\_ problem$

\eb

\end{document}